\newif\ifpersonal
\personaltrue 

\documentclass[a4paper, 10pt]{article}

\usepackage{preamble}
\usepackage[textwidth=20mm, textsize=tiny]{todonotes}
\usepackage{mathtools}
\usepackage{xr}
\usepackage[capitalize]{cleveref}

\title{A model for the assembly map of bordism-invariant functors}

\date{\today}

\author{Jordan Levin, Guglielmo Nocera, Victor Saunier}

\begin{document}
\maketitle

\begin{abstract}
    We study oplax colimits of stable categories, of hermitian categories and of Poincaré categories in nice cases. This allows us to produce a categorical model of the assembly map of a bordism-invariant functor of Poincaré categories which is also a Verdier projection, whose kernel we explicitly describe. 
    
    As a direct application, we generalize the Shaneson splitting for bordism-invariant functors of Poincaré categories proved by Calmès--Dotto--Harpaz--Hebestreit--Land--Moi--Nardin--Nikolaus--Steimle to allow for twists. We also show our methods can tackle their general twisted Shaneson splitting of Poincaré--Verdier localizing invariants which specifies to a twisted Bass--Heller--Swan decomposition for the underlying stable categories, generalizing part of recent work of Kirstein--Kremer.
\end{abstract}

\tableofcontents

\section{Introduction}

\hspace{1.2em} Let $E:\CatEx\to\Ecal$ be a localizing invariant of stable categories. Given a diagram $F$ in $\CatEx$ indexed by a space $X$, the assembly map of $E$ is the colimit-comparison map
$$
    \begin{tikzcd}
    \colim\limits_X E\circ F\arrow[r] &  E(\colim\limits_X F)
    \end{tikzcd}.
$$
If the source of $F$ is of the form $BG$, there are a number of conjectures, admitting partial positive results, which pertain to whether the assembly map is split-injective or an equivalence --- we refer to \cite{LuckIsoConjBook} for a very complete overview of those conjectures.

One of the reasons why the assembly map is hard to tract is that it is not \textit{a priori} the evaluation of a map in $\CatEx$. Using Efimov's breakthrough generalizing localizing invariants to dualizable categories \cite{EfimovDualizable}, Bartels--Efimov--Nikolaus \cite{BartelsEfimovNikolaus} have announced that the evaluation of any localizing invariant $E:\PrDual\to\Sp$ at a specific strongly-continuous functor
$$
    \begin{tikzcd}
        \Phi:\widehat{\mathrm{coShv}}(X; \Sp)\arrow[r] & \Sp^X
    \end{tikzcd}
$$
gives rise to the natural functor $E(\Phi):X\otimes E(\Sp)\to E(\Sp^X)$. The category $\widehat{\mathrm{coShv}}(X; \Sp)$ is however quite subtle to understand, being an internal mapping object in the category $\PrDual$. \\

The goal of this paper is to provide a different model of the assembly map, whose source admits a rather concrete description, in particular staying in the world of compactly-generated categories and such that the map in question is a Verdier projection, making investigations of its kernel possible. \\

However, the strategy does not quite work for localizing invariants out of $\CatEx$. Instead, we will show that given a diagram $F:X\to\CatP$ indexed by a space $X$ which is given a preferred choice of triangulation and valued in idempotent-complete Poincaré categories, there is a Poincaré--Verdier projection of Poincaré categories which models the assembly map for \textit{bordism-invariant} Verdier-localizing invariants of Poincaré categories. We refer to the series of papers \cite{HermKI, HermKII, HermKIII, HermKIV} for all what concerns Poincar\'e categories and hermitian K-theory. 

The main example of a bordism-invariant Verdier-localizing invariant is the L-theory functor $L:\CatP\to\Sp$, and its Karoubi-localizing version denoted either $\L$ or $L^{\langle -\infty\rangle}$. Bordism-invariance is critical to our strategy and being a purely hermitian notion, this approach is less enlightening for localizing invariants of stable categories; we also note that we do not need any enlargement to a Poincaré--dualizable categories setting, and this allows our strategy to work for Verdier-localizing invariants instead of Karoubi-localizing ones.

We also note that since quadratic functors of the form $\Qoppa^q_D$ are stable under colimits in $\CatP$, our model will also apply to the more familiar quadratic L-theory functor, which is given by restriction of the general L-theory functor to those Poincaré categories whose quadratic functor is of the form $\Qoppa^q_D$. There are more examples that can also arise from rings under suitable hypotheses of niceness, for instance the visible symmetric structure and, when it coincides with the latter, the symmetric structure. This in particular applies to the case of the $\tB\Z$-assembly or more generally discrete groups with no 2-torsion by Lemma 4.3.8 of \cite{HermKI}. \\

The central idea of our strategy is that localizing invariants are instances of decategorification. Hence, to understand the left hand side of the assembly map, it is worthwhile to investigate the 2-categorical version of the colimit. Such ideas have already been used with great success in \cite{LandTammePullbacks, LandTammePushouts} to understand pushouts under localizing invariants through the lax-pullback and the lax-pushout constructions. We will substantiate this approach through the study of \textit{oplax} colimits of Poincaré categories.

Our constructions find a precedent in the work of Ranicki and Ranicki--Weiss, albeit in a different language which was a precursor to the Poincaré categories that we use. More precisely in \cite[Proposition 13.7]{Ranicki-blue} and in Proposition 5.10 of \cite{RanickiWeiss}, they realize the left hand side of the $\L$-theoretic assembly map for additive 1-categories via a construction which is related to ours. Though we will not compare them precisely, we note that the latter is already featuring an idea that will play a crucial role for us, namely the use of barycentric subdivisions.

The oplax colimit of a functor $F:\Ical\to\Ccal$, where $\Ccal$ is a 2-category, is an object which satisfies a similar universal property of a colimit except cocones are replaced by \textit{oplax} cocones where now triangles are not required to commute but are only prescribed some natural transformation. Our convention of orientation is that oplax 2-simplices have their natural transformation $h\implies g\circ f$ going towards the composite. There is a similar notion of oplax limits satisfying dual conditions and we note that, when the indexing category is a space, there is no difference between a (co)limit and its oplax refinement. \\

In \cite{GepnerHaugsengNikolaus}, the authors identify the oplax colimits and the oplax limits in the 2-category $\cat{Cat}$ of categories. The latter is given by the category of sections of the cartesian fibration classifying the functor whose oplax limit we are considering in this paper, whereas the former is simply the total space of the cocartesian unstraightening.

It is a folklore result, of which we offer a proof, that given a strongly finite category (i.e. a finite category enriched in finite spaces), the oplax colimit and the oplax limit in the 2-category $\CatEx$ of stable categories coincide, and they are given by the oplax limit as in $\cat{Cat}$. In particular, such oplax (co)limits exist. We will also build in general the oplax colimit of $\CatEx$ and the map it receives from the cocartesian unstraightening. \\

The situation is significantly more complicated in $\CatP$, the 2-category of Poincaré categories. In the larger 2-category $\CatH$ of hermitian categories, it is rather straightforward to build an object having the correct universal property and such that its underlying stable category is the oplax colimit in $\CatEx$ of the underlying diagram --- we note however that the oplax colimit and the oplax limit no longer coincide in $\CatH$. We will show that if the diagram $F:\Ical\to\CatH$ is indexed by a \textit{strongly finite} category $\Ical$ and lands pointwise in hermitian categories whose quadratic functor is non-degenerate, then the quadratic functor on the hermitian oplax colimit $\UnH(F)$ is again non-degenerate.

It need not hold however that if $F$ lands in $\CatP$, then $\UnH(F)$ is a Poincaré category; we will show this fails already for the simplest of non-groupoid diagram shapes: $\Delta^1$. On the other hand, we give the following positive statement in Corollary \ref{UnPPoincareFiniteSimplicialComplex}, inspired from results of \cite{LurieNotesLtheory} and \cite{HermKI}:

\begin{prop}
    Let $K$ be a simplicial complex and denote $\Face(K)^\op$ the opposite of its poset of faces. Then, for every functor $F:\Face(K)^\op\to\CatP$, the hermitian oplax colimit $\UnH(F)$ is a Poincaré category.
\end{prop}

The proof of this result reduces to treating the case of the finite simplicial complexes associated to the $\Delta^n$, where one can provide a more concrete model to make the computation.

Still, note that even in this case, it will often not be the case that the oplax cocone giving its universal property to $\UnH(F)$ lands in $\CatP$, the failure being due to some functors $F(i)\to\UnH(F)$ not preserving the duality. \\

Given a simplicial complex $K$, its realization $|K|$ is equivalent to the space obtained from the opposite of its poset of faces $\Face(K)^\op$ by universally inverting all the arrows. The functoriality of oplax colimits therefore provides, for every $F:|K|\to\CatP$, a functor of hermitian categories
$$
    \begin{tikzcd}
        \Phi_F:\UnH(F\circ p)\arrow[r] & \colim\limits_{|K|} F
    \end{tikzcd}
$$
obtained by pulling along the localization $p:\Face(K)^\op\to|K|$. We prove the following in Corollary \ref{PVSequenceForOplaxColim}:

\begin{prop} \label{IntroPVSequence}
    Let $F:|K|\to\CatP$ where $K$ is a simplicial complex. There is a Poincaré--Verdier sequence of Poincaré categories:
    $$
        \begin{tikzcd}
            \ker(\Phi_F)\arrow[r] & \UnH(F\circ p)\arrow[r, "{\Phi_F}"] & \colim\limits_{|K|} F
        \end{tikzcd}
    $$
    In particular, $\ker(\Phi_F)$ is a Poincaré category.
\end{prop}

In particular, any Poincaré--Verdier-localizing invariant will induce a short-exact sequence whose third term is the right hand side of the assembly map. However, it need not be in general that $E(\UnH(F))$ is the left hand side of the assembly map; in fact, for non-necessarily bordism-invariant $E$, this value may depend on the choice of $K$ which realizes to $|K|$. \\

However, under the additional hypothesis of bordism-invariance we get in Theorem \ref{AssemblyMapOfBIPV} (see also Corollary \ref{AssemblyMapOfBIPVFinitary}):

\begin{thm}
    Let $F:|K|\to\CatP$ be valued in idempotent-complete Poincaré categories, where $K$ is a finite simplicial complex, and $E:\CatP\to\Ecal$ be a bordism-invariant Verdier-localizing functor. There is an exact sequence in $\Ecal$:
    $$
        \begin{tikzcd}
            E(\ker(\Phi_F))\arrow[r] & E(\UnH(F))\arrow[r, "{E(\Phi_F)}"] & E(\colim\limits_{|K|} F)
        \end{tikzcd}
    $$
    and an equivalence
    $$
        \begin{tikzcd}
            \colim\limits_{|K|} E\circ F\arrow[r, "{\simeq}"] & E(\UnH(F))
        \end{tikzcd}
    $$
    under which $E(\Phi_F)$ corresponds to the assembly map of $E$. If $E$ is further assumed to preserve filtered colimits, then the finiteness hypothesis on $K$ can be dropped.
\end{thm}

The key result to prove the above theorem is to show that for the standard simplices $\Delta^n$, the kernel of the Verdier sequence of Proposition \ref{IntroPVSequence} is actually a metabolic category; in particular, $E(\ker\Phi_F)=0$. \\

In Corollary \ref{AssemblyBZofBIPVisEquivalence}, we use the above to prove a version of the Shaneson splitting (originally \cite[Theorem 5.1]{Shaneson}) for bordism-invariant functors which allows for non-trivial coefficients --- this last part is precisely what we generalize from the Shaneson splitting found in \cite{HermKIV}. Note that though inspired by their strategy, our proof goes through a different simplicial complex than the one used in the proof of \cite{HermKIV}.

\begin{thm}[{Corollary \ref{AssemblyBZofBIPVisEquivalence}}]
    Let $E:\CatP\to\Ecal$ be a bordism-invariant Verdier-localizing functor and let $F:\tB\Z\to\CatP$ valued in idempotent-complete Poincaré-categories, then the assembly map
    $$
        \begin{tikzcd}
            \colim\limits_{\tB\Z} E\circ F\arrow[r, "{\simeq}"] & E(\colim\limits_{\tB\Z} F)
        \end{tikzcd}
    $$
    is an equivalence.
\end{thm}

This decomposition has a rich history, and can be seen as a refinement of the previously mentioned Shaneson splitting, but also, as in \cite{HermKIV}, a version of the Bass-Heller-Swan decomposition, originally found in K-theory in \cite[Theorem 2]{BHS-original} and \cite[The Fundamental Theorem (page 236)]{Grayson}, which is adapted to coefficients (i.e. suitably \textit{twisted}), as in \cite[§5, pages 427-428]{Ranicki-ESATS} and Grayson \cite[Theorem 2.3]{Grayson}. The most modern approach to the above statement in K-theory (and in fact, for localizing invariants of stable categories) is \cite{KirsteinKremer}.

As an application (Corollary \ref{Shaneson-Karoubi}), we obtain long exact sequence in Karoubi L-theory similar to Ranicki's long exact sequence for quadratic projective L-theory of twisted Laurent extensions \cite[Theorem 5.1]{Ranicki-ALTIII}.\\

 Working a bit harder, we are also able to show that the twisted Bass-Heller-Swan of \cite{HermKIV} can be transformed into a statement for Poincaré-localizing invariants, at least when the endomorphism is invertible. More precisely, we show in Theorem \ref{TSSVerdierLoc}:

\begin{thm}
    Let $F:\tB\Z\to\CatP$ be an idempotent-complete Poincaré category with $\Z$-action and $E:\CatP\to\Ecal$ a Poincaré--Verdier localizing invariant. Then, there is a split exact sequence
    $$
        \begin{tikzcd}
            \colim\limits_{\tB\Z}E\circ F\arrow[r] & E(\colim\limits_{\tB\Z} F)\arrow[r]\arrow[l, bend right=25, dotted] & NE^{\hyp}(F)
        \end{tikzcd}
    $$
    in $\Ecal$ whose first map is the assembly map of $E$ and whose second term, $NE^{\hyp}(F)$, only depends of the hyperbolisation of $E$ and in particular, vanishes if $E$ is bordism-invariant.
\end{thm}

The proof of the above theorem follows more closely the pattern of the non-twisted version of \cite{HermKIV}, and in particular we only recover the result of \cite{KirsteinKremer} for those localizing invariants of stable categories which factor through $\Hyp$. We also note that we make use of the fact that \cite{KirsteinKremer} show the decomposition is split for stable categories to deduce that there exists a splitting in the Poincaré setting, though we expect that trudging a little more, one can produce it directly.

\paragraph{Acknowledgements.}
We wish to thank Yonatan Harpaz for suggesting the problem and many helpful discussions throughout the making of this paper, as well as Christoph Winges for a careful reading and suggesting multiple improvements to a previous draft of this manuscript. We also thank the authors of \cite{HermKIV} for sharing an early draft of their manuscript.

We also warmly thank Arthur Bartels, Alexander Efimov, Dominik Kirstein, Christian Kremer, Dmitry Kubrak, Thomas Nikolaus, Fabian Hebestreit, Francesca Pratali, Phil Pützstück and Maxime Ramzi for fruitful discussions. This work was conducted under the auspices of the European Research Council as part of the project Foundations of Motivic Real K-Theory (ERC grant no. 949583). The second author was partially supported by the Simons Collaboration ``Perfection in Algebra, Geometry and Topology''.

\begin{conventions}
    Throughout this whole work, ``category'' will mean ``$\infty$-category'', and we adopt the framework of \cite{HTT, HA}. We also use the language of Poincaré categories as studied in \cite{HermKI, HermKII, HermKIII, HermKIV} and we refer the reader to those papers for the basic definitions and properties.

    We will assume that all our localizing invariants, whatever the type, are valued in a stable category.
\end{conventions}

\section{Oplax colimits of stable, hermitian and Poincaré categories}
\subsection{Recollections and oplax colimits for stable categories}

\hspace{1.2em} Let us begin a short recollection of general oplax colimits. We formulate the definitions in the context of \cite[Definition 6.1]{GepnerHaugsengNikolaus} --- this is weaker than usually required for $(\infty, 2)$-category theory, but we will not need more.

\begin{defi}
    An \textit{enhanced mapping functor} $\FunInt$ for a category $\Ccal$ is a functor
    $$
        \begin{tikzcd}
            \FunInt:\Ccal^\op\times\Ccal\arrow[r] & \Cat 
        \end{tikzcd}
    $$
    with a given natural equivalence of its underlying space $\FunInt^\simeq\simeq\Map_\Ccal$.
\end{defi}

\begin{ex} \label{ExampleEnhancedMappingFunctors}
    We will be mainly interested in the following examples:
    \begin{itemize}
        \item The category of functors $\Fun(\Ccal, \Dcal)$ determines an enhancement of the mapping spaces of $\Cat$.
        \item The category of exact functors $\FunEx(\Ccal, \Dcal)$ determines an enhancement of the mapping spaces of $\CatEx$.
        \item The category of hermitian functors $\FunH((\Ccal,\Qoppa), (\Dcal,\Psi))$, the total space of the left fibration classifying $\FunEx(\Ccal, \Dcal)\to\Spaces$ sending $f$ to $\Nat(\Qoppa, \Psi\circ f)$, determines an enhancement of the mapping spaces of $\CatH$.
        \item The full subcategory $\FunP((\Ccal,\Qoppa), (\Dcal,\Psi))$ of $\FunH((\Ccal,\Qoppa), (\Dcal,\Psi))$ spanned by those hermitian functors $(f, \eta)$ which are Poincaré determines an enhancement of the mapping spaces of $\CatP$.
    \end{itemize}
    The last category is less well-behaved than the other three, see \cite[Section 6.5]{HermKI}.
\end{ex}

\begin{defi} \label{DefOplaxColimit}
    Let $(\Ccal, \FunInt)$ be a category with a choice of enhanced mapping functor. An oplax cocone of a functor $F:\Ical\to\Ccal$ is a point in the category 
    $$
        \Fun_{\Ical}(\Ical, \UnCart(\FunInt(F(-), X)))
    $$
    for some $X\in\Ccal$. By varying $X$, the above category determines a functor $\Ccal\to\CatInfty$ whose unstraightening is the category of oplax cocones of $F$. \\
    
    The \textit{oplax colimit} of a functor $F$, if it exists, is an oplax cocone which is initial in the enhanced sense, so that there is an equivalence
    $$
        \begin{tikzcd}
            \FunInt(\oplaxcolim(F), X)\arrow[r, "\simeq"] & \Fun_{\Ical}(\Ical, \UnCart(\FunInt(F(-), X)))
        \end{tikzcd}
    $$
    of categories, induced by post-composing the oplax cocone of $\oplaxcolim(F)$.
\end{defi}

We will be interested in the interaction of oplax colimits and precomposition of the diagram. We have the following lemma:

\begin{lmm} \label{OplaxColimitsPreserveColimits}
    Let  $(\Ccal, \FunInt)$ be a category with a choice of enhanced mapping functor, such that $\Ccal$ is cocomplete and $(\Ccal, \FunInt)$ has all (small) oplax colimits. Then, given a functor $F:K\to\Ccal$, the association
    $$
        \alpha\in\CatInfty_{/K}\longmapsto \oplaxcolim(F\circ\alpha)\in\Ccal
    $$
    is functorial and preserves colimits.
\end{lmm}
\begin{proof}
    Since we have assumed the existence of both colimits and oplax colimits, it suffices to work with represented functors. In particular, the Yoneda lemma reduces both the functoriality and the preservation of colimits to checking that, for a fixed functor $F:\Ccal\to\CatInfty$, the association
    $$
        (\alpha:\Dcal\to\Ccal)\mapsto \Fun_\Dcal(\Dcal, \Un(F\circ\alpha))\simeq\Fun_{/\Ccal}(\Dcal, \Un(F))
    $$
    is functorial and sends colimits to limits, which is straightforward from identifying the right hand side with the mapping category from $\alpha$ to $(\Un(F)\to\Ccal)$ in $\CatInfty_{/\Ccal}$. Here, note that we need only appeal to the 1-categorical Yoneda since we are only proving statements related to 1-colimits and 1-functoriality.
\end{proof}

 Given a functor $F:\Ical\to\CatInfty$, there are two fibrations associated to it via the unstraightening. We write $\Un(F)$ for the total space of its cocartesian unstraightening and $\UnCart(F)$ for the cartesian one.

\begin{prop} \label{FunctorsFromUn}
    Let $F:\Ical\to\CatInfty$ be a functor. For every category $\Ecal$, there are equivalences:
    $$
    \begin{tikzcd}[row sep=0.2em]
        \Fun(\Un(F), \Ecal)\arrow[r, "{\simeq}"] & \Fun_{\Ical}\left(\Ical, \UnCart(\Fun(F(-), \Ecal))\right) \\
        \Fun(\Un(F)^\op, \Ecal)\arrow[r, "{\simeq}"] & \Fun_{\catop{\Ical}}\left(\catop{\Ical}, \Un(\Fun(\catop{F(-)}, \Ecal))\right)
    \end{tikzcd}
    $$
    where, on the right, $\Fun_{\Ical}(\Ical, -)$ denotes the category of sections of a given fibration over $\Ical$. 
    
    In particular, the collection of inclusions of the fibers $F(i)\to\Un(F)$ upgrades to an oplax cocone which realizes $\Un(F)$ as the oplax colimit of $F:\Ical\to\CatInfty$.
\end{prop}
\begin{proof}
    By \cite[Corollary 3.2.2.13]{HTT} and Proposition 7.3 of \cite{GepnerHaugsengNikolaus}, the functor $\Fun(F(-), \Ecal):\catop{\Ical}\to\CatInfty$ is classified by a cartesian fibration $\UnCart(\Fun(F(-), \Ecal))\to\Ical$ which satisfies the following equivalence
    $$
        \begin{tikzcd}
            \Fun(\Un(F), \Ecal)\arrow[r, "{\simeq}"]& \Fun_{\Ical}\left(\Ical, \UnCart(\Fun(F(-), \Ecal))\right)
        \end{tikzcd}
    $$
    which is induced by the inclusions of the fibers. Let us also explain a model-independent argument: by the equivalence of $(\infty,2)$-categories of straightening--unstraightening \cite[Theorem 3.2.0.1]{HTT}, we have
    $$
        \Fun_{\Ical}\left(\Ical, \UnCart(\Fun(F(-), \Ecal))\right)\simeq\Nat^{\oplax}(\cst_*, \Fun(F(-), \Ecal))
    $$
    where $\Nat^{\oplax}$ is the category of oplax natural transformations, i.e. the category of maps in the $(\infty, 2)$-category of oplax functors $\Ical\to\Cat$ or, equivalently, in that of cartesian fibrations over $\Ical$ (see section 5 of \cite{HaugsengHebestreitLinskensNuiten} for details). Of course 
    $$
        \Nat^{\oplax}(\mathrm{cst}_*, \Fun(F(-), \Ecal))\simeq\Nat^{\oplax}(F(-), \Ecal)
    $$ 
    and the covariant unstraightening gives
    $$
        \Nat^{\oplax}(F(-), \Ecal)\simeq\Fun_{\Ical}\left(\Un(F), \Ecal\times\Ical\right)\simeq\Fun(\Un(F), \Ecal).
    $$
    The second equivalence can be deduced from the first by using that the opposite of a cartesian fibration classifying $F$ is the cocartesian fibration classifying $\mathrm{op}\circ F$ (see for instance \cite{BarwickGlasmanNardin}), so that we have:
    $$
        \begin{tikzcd}
            \Fun(\Un(F)^\op, \Ecal)\arrow[r, "{\simeq}"] & \Fun_{\catop{\Ical}}\left(\catop{\Ical}, \Un(\Fun(\catop{F(-)}, \Ecal))\right)
        \end{tikzcd}
    $$
    which concludes the proof.
\end{proof}

\begin{rmq} \label{ComparisonWithGHN}
    In \cite{GepnerHaugsengNikolaus}, a different definition of oplax colimits (as well other variants of (op)lax (co)limits) is given through weighted limits. It follows from the results of section 7 of \loccit but also from the above proposition that these agree with our definition.
    
    Note also that \cite[Proposition 7.1]{GepnerHaugsengNikolaus} identifies the oplax limit of a functor $F:\Ical\to\CatInfty$ with the category of sections $\Fun_{\Ical^\op}(\Ical^\op, \UnCart(F))$. Indeed, by a similar proof as in proposition \ref{FunctorsFromUn}, there is an equivalence
    $$
        \Fun(\Ecal, \Fun_{\Ical^\op}(\Ical^\op, \UnCart(F)))\simeq\Fun_{\Ical^\op}(\Ical^\op, \UnCart(\Fun(\Ecal, F(-))).
    $$
    Therefore, the first equivalence of Proposition \ref{FunctorsFromUn} can be restated as the following:
    $$
        \begin{tikzcd}
            \Fun(\oplaxcolim(F), \Ecal)\arrow[r, "{\simeq}"] & \oplaxlim(\Fun(F(-), \Ecal))
        \end{tikzcd}
    $$
    which is an oplax-colimit variant of the fact that the contravariant entry of mapping spaces sends colimits to limits.
\end{rmq}

Suppose now that our functor $F$ lands in $\CatEx$, the category of stable categories and exact functors between them. In general, $\Un(F)$ will not be a stable category itself --- this already fails if $F$ is constant and $\Ical$ not stable. We now seek to produce a stable category as close as possible to $\Un(F)$.

\begin{defi}
    Let $\Ecal$ be a stable category, and $F:\Ical\to\CatEx$ a functor. A functor $G:\Un(F)\to\Ecal$ is called \textit{fiberwise-exact} if its restriction along each fiber $G_i:F(i)\to\Un(F)\rightarrow\Ecal$ is an exact functor. We denote by
    $$
        \Fun^{\fbwex}(\Un(F), \Ecal)
    $$
    the full subcategory of $\Fun(\Un(F), \Ecal)$ of fiberwise-exact functors.
\end{defi}

\begin{thm} \label{StabUnDefThm}
    Let $\Ical$ be a category and $F:\Ical\to\CatEx$ a functor. There exists a stable category $\StabUn(F)$ equipped with a fiberwise-exact functor $\alpha:\Un(F)\to\StabUn(F)$ satisfying the following universal property: for every stable $\Ecal$, the functor
    $$
        \begin{tikzcd}
            \alpha^*:\FunEx(\StabUn(F), \Ecal)\arrow[r, "\simeq"] & \Fun^{\fbwex}(\Un(F), \Ecal)
        \end{tikzcd}
    $$
    is an equivalence.
\end{thm}
\begin{proof}
    This follows from the general framework of \cite[Proposition 5.3.6.2]{HTT} as we now explain. Let $\Pcal(\Un(F))$ be the category of space-valued presheaves on $\Un(F)$. Consider the full subcategory $\Pcal^{\fbwex}(\Un(F))$ of $\Pcal(\Un(F))$ spanned by those presheaves whose restriction to each fiber is right-exact; this is an accessible localization of $\Pcal(\Un(F))$ and we denote by $L$ the associated localization functor. 
    
    Let $j$ be the Yoneda embedding $\Un(F)\to\Pcal(\Un(F))$ and write $\Qcal$ for the smallest full subcategory of $\Pcal^{\fbwex}(\Un(F))$ stable under finite colimits and containing the image of $L\circ j$. It follows from \cite[Proposition 5.3.6.2]{HTT} that, for every category with finite colimits $\Ecal$, the functor $L\circ j:\Un(F)\to\Qcal$ induces an equivalence
    $$
        \begin{tikzcd}
            (L\circ j)^*:\FunLEx(\Qcal, \Ecal)\arrow[r, "\simeq"]& \Fun^{\fbwex}(\Un(F), \Ecal)
        \end{tikzcd}
    $$
    with source the category of left-exact functors $\Qcal\to\Ecal$. Thus, to conclude we can simply let $\StabUn(F)$ be the Spanier-Whitehead stabilization of $\Qcal$, that is the image of $\Qcal$ under the left adjoint of the inclusion $\CatEx\subset\cat{Cat}^{\mathrm{fin-colim}}$.
\end{proof}

\begin{prop} \label{StabUnIsOplaxColimit}
    Let $F:\Ical\to\CatEx$ be a functor. Then, $\StabUn(F)$ is the oplax colimit of $F$ computed in $\CatEx$ with its enhanced mapping functor as in Example \ref{ExampleEnhancedMappingFunctors}.
\end{prop}
\begin{proof}
    This follows without difficulty from the previous results. Let $\Ecal$ be stable, then Theorem \ref{StabUnDefThm} shows that the map $\Un(F)\to\StabUn(F)$ induces 
    \[
        \begin{tikzcd}
            \FunEx(\StabUn(F), \Ecal)\arrow[r, "\simeq"]&  \Fun^{\fbwex}(\Un(F),\Ecal)
        \end{tikzcd}.
    \]
    
    Since being fiberwise-exact is a fiberwise condition, as the name suggests, we deduce from Theorem \ref{FunctorsFromUn} that pulling along the oplax cocone of $\Un(F)$ induces
    \[
        \begin{tikzcd}
            \Fun^{\fbwex}(\Un(F),\Ecal) \arrow[r, "\simeq"] & \Fun_{\Ical}(\Ical, \UnCart(\Fun^{\Ex}(F(-), \Ecal)))
        \end{tikzcd}
    \]
    which concludes the proof.
\end{proof}

\begin{ex}
    If $F:\Ical\to\CatEx$ is constant with value $\Ccal$, then $\Un(F)\simeq\Ical\times\Ccal$ and $\StabUn(F)$ is $\Ical\otimes\Ccal$, the value of the left adjoint of the endofunctor of $\CatEx$ given by $\Fun(\Ical, -)$: this is often called the tensor of the stable category $\Ccal$ by $\Ical$ in the sense of $\CatEx$, and was described under the name \textit{semi-exact tensoring} in \cite[Definition 2.1]{SaunierFundamental}.
\end{ex}

\begin{ex}
    If $F:\Delta^1\to\CatEx$ is a functor corresponding to an exact functor $G:\Bcal\to\Acal$, then $\StabUn(F)$ is the total space of the semi-orthogonal decomposition classified by $G$, i.e. it fits in the pullback square:
    $$
        \begin{tikzcd}
            \StabUn(F)\arrow[d]\arrow[r] & \Acal^{\Delta^1}\arrow[d, "t"] \\
            \Bcal\arrow[r] & \Acal
        \end{tikzcd}.
    $$
    Indeed, this follows from \cite[Remark 5.7]{ChristDycherhoffWaldeLaxAdd} and the identification of $\StabUn(F)$ with the oplax colimit of $F$.
\end{ex}

\begin{rmq} \label{InverseToUniversalPropEquivalence}
    In Lemma \ref{AlphaiasAdjoints} and Proposition \ref{ExplicitStabUnStronglyFinite}, we will show that if $\Ical$ is strongly finite (see Definition \ref{DefStronlyFinite}), the inverse to the first equivalence of Proposition \ref{FunctorsFromUn} is given by right Kan extension along $\alpha$. 
\end{rmq}

Let us give a more explicit description of $\StabUn(F)$ by using Proposition \ref{FunctorsFromUn}, although we will not use it in the following. It follows from the Proposition that $\Pcal(\Un(F))$ is equivalent to the category of sections of $\Un(\Fun(\catop{F(-)}, \Spaces))$ and the reflexive subcategory $\Pcal^{\fbwex}(\Un(F))$ corresponds exactly to the following category of sections:
$$
    \Pcal^{\fbwex}(\Un(F))\simeq\Fun_{\catop{\Ical}}(\catop{\Ical}, \Un(\FunLEx(\catop{F(-)}, \Spaces))).
$$
We can write the right-hand side more gently as
$$
    \Fun_{\catop{\Ical}}(\catop{\Ical}, \Un(\Ind(F)))
$$

where $\Ind(-)$ has the left-Kan-extension functoriality (hence covariant).
This perspective shows that $\Pcal^{\fbwex}(\Un(F))$ is already stable, by combining \cite[Definition 2.9 and Corollary 7.7]{GepnerHaugsengNikolaus}. In particular, we get that $\StabUn(F)$ can be identified with the smallest stable subcategory of $\Pcal^{\fbwex}(\Un(F))$ containing the image of $\Un(F)$ under the reflexion $\Pcal(\Un(F))\to\Pcal^{\fbwex}(\Un(F))$.

\subsection{Oplax colimits of hermitian categories}

We will be interested in functors valued in the category of Poincaré categories, $\CatP$. Via the forgetful functor $\fgt:\CatP\to\CatEx$, all such functors have in particular an underlying functor to $\CatEx$ to which the previous section applies. We now seek to upgrade the structure on $\StabUn(F)$ accordingly; let us first deal with an intermediate case: namely, a functor $F:\Ical\to\CatH$, where $\CatH$ denotes the category of hermitian categories. \\

Fix $F:\Ical\to\CatEx$. We are interested in how much datum is a dotted lift
$$
    \begin{tikzcd}
        & \CatH\arrow[d] \\
        \Ical\arrow[r, "F"]\arrow[ru, dotted] & \CatEx
    \end{tikzcd}.
$$
The following is straightforward:

\begin{lmm} \label{DatumOfLift}
    There is an equivalence of categories
    $$
        \Fun_{/\CatEx}(\Ical, \CatH)\simeq\Fun_{\Ical}(\Ical, \UnCart(\Fun^q(F(-)))).
    $$
\end{lmm}
\begin{proof}
    By definition, $\CatH\simeq\UnCart(\Fun^q(-))$ so that $\UnCart(\Fun^q(F(-)))\to\Ical$ is the pullback of the cartesian fibration $\CatH\to\CatEx$ along $F$. This concludes.
\end{proof}

In particular, such a lift of $F$ gives quadratic functors $\Qoppa_{F(i)}:\catop{F(i)}\to\Sp$ for $i\in\Ical$ and, for each $f:i\to j$ in $\Ical$, arrows $\eta_f:\Qoppa_{F(i)}\to\Qoppa_{F(j)}$ in the cartesian unstraightening. The latter are equivalently given by natural transformations of quadratic functors as follows
$$
    \begin{tikzcd}
        \Qoppa_{F(i)}\arrow[r] & \Qoppa_{F(j)}\circ F(f)^\op
    \end{tikzcd}
$$
or, equivalently, $F(f)^\op_!\Qoppa_{F(i)}\to\Qoppa_{F(j)}$. \\

Given $F:\Ical\to\CatEx$, we write $\alpha_i$ for the exact functors obtained as the composites 
$$
    \begin{tikzcd}
        \alpha_i:F(i)\arrow[r] & \Un(F)\arrow[r] & \StabUn(F)
    \end{tikzcd}
$$
which form together the oplax cocone of $\StabUn(F)$. Given a lift of $F$ as above and $i\in\Ical$, there is a quadratic functor on $\StabUn(F)$ obtained by left Kan extending $\Qoppa_{F(i)}$ along $\alpha_i^\op$; this association is functorial in $\Ical$ as we now explain.

For an arrow $f:i\to j$ in $\Ical$, we have a natural transformation $\eta_F(f):(F(f)^\op)_!\Qoppa_{F(i)}\to\Qoppa_{F(j)}$ supplied by the lift of $F$, as well as a natural transformation $\alpha_i\to\alpha_j\circ F(f)$ induced by the oplax cocone defining $\StabUn(F)$.


Combining them, one gets a natural transformation of functor $\StabUn(F)^\op \to \Sp$
$$
    \begin{tikzcd}[column sep=large]
        (\alpha_{i}^\op)_!\Qoppa_{F(i)}\arrow[r] & (\alpha_{j}^\op\circ F(f)^\op)_!\Qoppa_{F(i)}\arrow[r, "{(\alpha_{j}^\op)_!\eta_F(f)}"] & (\alpha_{j}^\op)_!\Qoppa_{F(j)}
    \end{tikzcd}
$$
which is how we want the functoriality to look on arrows. 
Let us now work more formally:

\begin{lmm}\label{ExistenceOfQ}
    Suppose $\widetilde F:\Ical\to\CatH$ is a functor lifting $F:\Ical\to\CatEx$. Then, there is a functor $\Qcal:\Ical\longrightarrow\Fun^q(\StabUn(F))$ given on objects by
    $$
        i\in\Ical\longmapsto(\alpha_i^\op)_!\Qoppa_{F(i)}
    $$
    and on arrows by the maps described previously.
\end{lmm}
\begin{proof}
    We will do more than the statement asks and prove that $\Qcal$ can be obtained functorially from the lift $\widetilde F$, namely as the value of a functor
    $$
        \begin{tikzcd}
            \Theta:\Fun_{\Ical}(\Ical, \UnCart(\Fun^q(F(-))))\arrow[r] & \Fun(\Ical,\Fun(\StabUn(F)^\op, \Sp))
        \end{tikzcd}
    $$
    Note that since $(\alpha_i^\op)_!\Qoppa_{F(i)}$ is automatically quadratic by \cite[Lemma 1.4.1 (iii)]{HermKI} and quadratic functors form a \textit{full} subcategory, we may omit the superscript $q$ from the target of $\Theta$.
    
    Now, observe that, since all the $\alpha_i:F(i)\to\StabUn(F)$ factor through the common $\alpha:\Un(F)\to\StabUn(F)$ and left Kan extending is a functorial operation, it suffices to produce a functor:
    $$
        \begin{tikzcd}
            \Phi:\Fun_{\Ical}(\Ical, \UnCart(\Fun^q(F(-))))\arrow[r] & \Fun(\Ical,\Fun(\Un(F)^\op, \Sp)).
        \end{tikzcd}
    $$
    which takes a section $i \mapsto \Qoppa_{F(i)}$ and maps it to the functor $i\mapsto(\beta_i^\op)_!\Qoppa_{F(i)}$ where $\beta_i:F(i)\to\Un(F)$ is the inclusion of the fiber. Indeed, $\widetilde F$ defines an object on the left hand side coming by Lemma \ref{DatumOfLift}. \\

    We will show that $\Phi$ is actually given by postcomposition along some $\Psi$, where $\Psi$ is a functor
    $$
        \begin{tikzcd}
            \Psi:\UnCart(\Fun^q(F(-)))\arrow[r] & \Fun(\Un(F)^\op, \Sp)
        \end{tikzcd}
    $$
    mapping $\Qoppa$ to $(\beta_i^\op)_!\Qoppa$ where $i$ is the point in the fiber over which $\Qoppa$ lives. Finally, we remark that
    $$
        \UnCart(\Fun^q(F(-)))\simeq\Un(\Fun^q(F(-)))
    $$
    where functoriality on the right hand side is now given by left Kan extension. Using Proposition \ref{FunctorsFromUn}, supplying $\Psi$ reduces to explain that there is a point in
    $$
        \Fun_{\Ical}(\Ical, \UnCart(\Fun(\Fun^q(F(-)), \Fun(\Un(F)^\op,\Sp)))))
    $$
    which sends $i$ to the object of $\Fun(\Fun^q(F(i)), \Fun(\Un(F)^\op,\Sp)))$ which takes a quadratic functor $\Qoppa:\catop{F(i)}\to\Sp$ and left Kan extends along $\beta_i^\op$. We remark that since 
    $$
        \id_{\Un(F)}\in\Fun(\Un(F), \Un(F))\simeq\Fun_{\Ical}(\Ical, \UnCart(\Fun(F(-),\Un(F))))
    $$
    and $\Fun_{\Ical}(\Ical, -)$ is functorial with source $\Cat_{/\Ical}$, it suffices to show that there is a dashed map 
    $$
        \begin{tikzcd}
             \UnCart(\Fun(F(-),\Un(F)))\arrow[rd]\arrow[rr, dotted] & & \UnCart(\Fun(\Fun^q(F(-)), \Fun(\Un(F)^\op,\Sp)))\arrow[ld] \\
                &\Ical
        \end{tikzcd}
    $$
    making the diagram commute and sending $\beta:F(i)\to\Un(F)$ to the left Kan extension functor $\beta_!:\Fun^q(F(i))\to\Fun(\Un(F)^\op,\Sp)$. Note that this map will not preserve cocartesian edges: this is precisely why we picked up the second arrow of the composite in the functoriality we sketched before the statement of the Lemma. Nevertheless, it suffices to produce this map after unstraightening, so that we have to provide
    $$
        \begin{tikzcd}
            \Fun(F(i),\Un(F))\arrow[r] & \Fun(\Fun^q(F(i)), \Fun(\Un(F)^\op,\Sp))
        \end{tikzcd}
    $$
    of functors natural in $i\in\Ical$, which is given pointwise by the functor sending $f$ to $f_!$.  
    
    At this point in the proof, we observe that we have arrived at a simple statement: the association $f\mapsto f_!$ is functorial in both $f$ and the source category. Indeed, it is straightforward that $f\mapsto f^*$ is functorial in this way, where $f^*$ denotes the precomposition functor; consequently, its left adjoint also shares its functoriality by a direct application of the Yoneda lemma. This functoriality is precisely the first arrow in the composite we have written. This concludes the proof of Lemma \ref{ExistenceOfQ}.
\end{proof}

\begin{defi} \label{definition-Qoppa-Ex}
    Suppose $F:\Ical\to\CatH$ is a functor. Then we write $\Qoppa_F$ for the quadratic functor:
    $$
        \Qoppa_F\coloneqq\left(\colim_{i\in\Ical}\Qcal(i)\right)\in\Fun^q(\StabUn(F))
    $$
    We denote $\UnP(F)$ the hermitian category $(\StabUn(F), \Qoppa_F)$.
\end{defi}

\begin{prop}
    Let $\Qoppa'$ be a quadratic functor on $\StabUn(F)$, then there is a natural equivalence as follows:
    $$
        \begin{tikzcd}
        \Nat(\Qoppa_F, \Qoppa')\arrow[r, "{\simeq}"] & \Fun_{\Ical}(\Ical, \UnCart(\Nat(\Qcal(-), \Qoppa'))).
        \end{tikzcd}
    $$
    In particular, the left-hand-side is a space, and therefore $\Qoppa_F$ is not only the colimit (by definition) but also the oplax colimit of $\Qcal$ (for trivial reasons).
\end{prop}
\begin{proof}
    By \cite[Corollary 3.3.3.4]{HTT}, since natural transformations form a space, there is an equivalence
    $$
        \begin{tikzcd}
        \Fun_{\Ical}(\Ical, \UnCart(\Nat(\Qcal(-), \Qoppa')))\arrow[r, "{\simeq}"] & \displaystyle\lim_{i\in\Ical}\Nat(\Qcal(i), \Qoppa')
        \end{tikzcd}.
    $$
    In particular, it is clear that $\Qoppa_F$ has the wanted universal property with respect to quadratic functors on $\StabUn(F)$.
\end{proof}

\begin{rmq} \label{IndDualityAsSection}
    It follows from the above that the associated bilinear $B_{\Qoppa_F}$ and the $\Ind$-duality $D_{\Qoppa_F}$ --- i.e. the functor $\StabUn(F)^\op\to\Ind(\StabUn(F))$ formally associated to $B_{\Qoppa_F}$ via currying --- also have similar universal properties. We also get that they are given by the formulas:
    \begin{align*}
        B_{\Qoppa_F}\coloneqq\colim_{i\in\Ical}(\alpha_i^\op\times\alpha_i^\op)_!B_{\Qoppa_{F(i)}} & & 
        D_{\Qoppa_F}\coloneqq\colim_{i\in\Ical}(\alpha_i^\op)_!\left(\Ind(\alpha_i)\circ D_{\Qoppa_{F(i)}}\right)
    \end{align*}
    where the colimits are taken in the correct functor categories, and consequently also pointwise.
\end{rmq}

\begin{thm}
    Let $F:\Ical\to\CatH$ be a functor. Then, the hermitian category $\UnP(F)$ is the oplax colimit of $F$ in the category $\CatH$ with the enhanced mapping functor of Example \ref{ExampleEnhancedMappingFunctors}, in the sense of Definition \ref{DefOplaxColimit}.
\end{thm}
\begin{proof}
    We have to show that for every hermitian category $(\Ccal, \Qoppa)$, there is a (natural) equivalence:
    $$
        \begin{tikzcd}
            \FunH(\UnP(F), (\Ccal, \Qoppa))\arrow[r, "{\simeq}"] & \Fun_{\Ical}(\Ical, \UnCart(\FunH(F(-), (\Ccal, \Qoppa))))
        \end{tikzcd}.
    $$
    Recall that for every pair of hermitian categories $(\Ccal, \Qoppa), (\Dcal, \Qoppa')$, there is a functor 
    $$
        \begin{tikzcd}
            \FunH((\Ccal, \Qoppa),(\Dcal, \Qoppa'))\arrow[r] & \FunEx(\Ccal, \Dcal)  
        \end{tikzcd}
    $$
    which is a left fibration whose fiber at a functor $f:\Ccal\to\Dcal$ is given by $\Nat(\Qoppa, \Qoppa'\circ\catop{f})$. In particular, for every exact $f:\StabUn(F)\to\Ccal$, we get a diagram whose columns are fiber sequences as follows
    $$
        \begin{tikzcd}
            \Nat(\Qoppa_F, \Qoppa\circ\catop{f})\arrow[r]\arrow[d] & \Fun_{\Ical}(\Ical, \UnCart(\Nat(\Qoppa_{F(-)}, \Qoppa\circ\catop{f}\circ\alpha_{(-)}^\op)))\arrow[d] \\ 
            \FunH(\UnP(F), (\Ccal, \Qoppa))\arrow[r]\arrow[d] & \Fun_{\Ical}(\Ical, \UnCart(\FunH(F(-), (\Ccal, \Qoppa))))\arrow[d] \\
            \FunEx(\StabUn(F), \Ccal)\arrow[r] & \Fun_{\Ical}(\Ical, \UnCart(\FunEx(F(-), \Ccal)))
        \end{tikzcd}
    $$
    where $\alpha_{(-)}$ denotes the functor $F(-)\to\StabUn(F)$. By Proposition \ref{StabUnIsOplaxColimit}, the bottom horizontal arrow is an equivalence, and since the bottom vertical maps are left fibrations, it suffices to show that the top arrow is an equivalence for every $f$. This latter fact has been proven in the previous Proposition (apply to $\Qoppa'=\Qoppa\circ\catop{f}$), so this concludes the proof.
\end{proof}

\subsection{Oplax colimits of Poincaré categories}

\hspace{1.2em} Suppose now given $F:\Ical\to\CatP$; it need not be the case that the hermitian category $\UnP(F)$ is Poincaré. Already for constant $F$, this does not hold for general $\Ical$, see \cite[Sections 6.4, 6.5]{HermKI}. Let us try to explain why this is to be expected. Recall that $\CatP$ fits in a pullback square
$$
    \begin{tikzcd}
        \CatP\arrow[r]\arrow[d] & (\CatH)_{\mathrm{non-deg}}\arrow[d] \\
        (\CatEx)^{\htpyfixed\cyclicGrp{2}}\arrow[r] & (\CatEx)^{\mathrm{lax}-\htpyfixed\cyclicGrp{2}}
    \end{tikzcd}
$$
where the $\cyclicGrp{2}$-action on $\CatEx$ is given by taking opposites. The bottom left corner is defined to be $\cyclicGrp{2}$-homotopy fixed points for this action. Concretely, an object of this category is a \textit{stable category with duality}, a stable category $\Ccal$ with an equivalence $\begin{tikzcd}[cramped]D:\Ccal^\op\arrow[r, "{\simeq}"] & \Ccal\end{tikzcd}$ along with coherence data, and maps between them being exact functors intertwining such equivalences. We may view the left vertical arrow as the factor which takes a Poincaré category and sends to the associated stable category with duality. 

The term $(\CatEx)^{\mathrm{lax}-\htpyfixed\cyclicGrp{2}}$ is a bit more mysterious: by this we mean the category whose objects are stable categories $\Ccal$ equipped with merely an exact functor $D:\Ccal^\op\to\Ccal$, and the morphisms are given by exact functors $\Ccal\to\Dcal$ and a non-necessarily invertible natural transformation $D_{\Dcal}\circ\catop{f}\to f\circ D_\Ccal$; the right vertical functor is thus passing to the duality of a non-degenerate hermitian category.

The subtleties concerning oplax colimits in $\CatP$ stem from the bottom left corner: indeed, the $\cyclicGrp{2}$-homotopy fixed points is not a 2-categorical construction in that it imposes the restriction of \textit{invertible} natural transformations. In contrast, the right hand side of the square involves non-invertible natural transformations by its very definition. In this section, we will see that under suitable finiteness hypotheses, guaranteeing non-degeneracy of the duality will be rather straightforward, but we will have to impose stringent conditions to get perfectness of the duality and duality-preservation of the induced maps between oplax colimits. \\

The following definition, lifted verbatim from Definition 6.5.1 of \cite{HermKI}, captures a sufficient condition on $\Ical$ for the oplax colimits of $(\CatH)_{\mathrm{non-deg}}$ to be computed as in $\CatH$:

\begin{defi} \label{DefStronlyFinite}
    A category $\Ical$ is \textit{strongly finite} if it is finite and for all $i,j\in\Ical$, the mapping space $\Map_\Ical(i, j)$ is also finite.
\end{defi}

For the purpose of indexing (op)lax colimits, strongly finite categories are the higher categorical equivalents of finite spaces; and indeed, they are the finite objects of $\CatInfty$ which are further enriched in finite spaces, hence a reasonable categorification. 

Given a functor $F:\Ical\to\CatEx$, we begin by investigating $\Fun_{\catop{\Ical}}(\catop{\Ical}, \UnCart(F))$. By \cite[Proposition 7.1]{GepnerHaugsengNikolaus} this is the oplax limit of the associated functor $I \to \CatEx \to \Cat$. 

Given a section $\phi$, one can check that the section $\Omega\phi$ is the section given pointwise by $\Omega_{F(j)}\phi(j)$ and, since all the $F(j)$ are stable, $\Omega$ is invertible with inverse given by taking a pointwise suspension in a similar manner. Hence $\Fun_{\catop{\Ical}}(\catop{\Ical}, \UnCart(F))$ is already stable. We will show that when $\Ical$ is strongly finite, $\Fun_{\catop{\Ical}}(\catop{\Ical}, \UnCart(F))$ is a also model for the oplax colimit in $\CatEx$ of the functor $F$.

\begin{lmm}\label{AlphaiasAdjoints}
    Suppose that all of the mapping spaces of $\Ical$ are finite. Then, for all $i\in\Ical$, the functor $\ev_i:\Fun_{\catop{\Ical}}(\catop{\Ical}, \UnCart(F))\to F(i)$ admits a right adjoint, given in formula, for $X\in F(i)$: 
    \begin{align*}
        \alpha_i(X):j\in\Ical^\op\longmapsto \left(\lim_{\gamma\in\Map_{\Ical}(i, j)}F(\gamma)(X)\right)\in F(j)\subset\UnCart(F).
    \end{align*}
    Moreover, via Proposition \ref{FunctorsFromUn}, the collection of right adjoints upgrades to a functor
    $$
        \alpha:\Un(F)\to\Fun_{\catop{\Ical}}(\catop{\Ical}, \UnCart(F))
    $$ 
    whose restrictions along each fiber are precisely given by $\alpha_i$.
\end{lmm}
\begin{proof}
    We begin with the first statement. For brevity, we will write $\Gamma$ for $\Fun_{\catop{\Ical}}(\catop{\Ical}, \UnCart(F))$ in the rest of the proof. Identifying $F(i)\simeq\Fun_{/\Ical^\op}(*, \UnCart(F))$ where $*$ is viewed over $\Ical^\op$ via the inclusion $\{i\}\subset\Ical^\op$, $\ev_i$ is given by precomposition along $i:*\to\Ical^\op$. Hence, if it exists, the wanted right adjoint must be given by the relative right Kan extension in the diagram:
    $$
        \begin{tikzcd}
            *\arrow[r,"X"]\arrow[d, "i"'] & \UnCart(F)\arrow[d] \\
            \Ical^\op\arrow[r, equal]\arrow[ru, dotted] & \Ical^\op
        \end{tikzcd}.
    $$
    The proposed formula is the putative formula for such relative Kan extensions which can be seen from combining (the dual statement to) \cite[Corollary 4.3.2.14]{HTT} and the fact that limits in $\UnCart(F)$ can be computed by pushing along the cartesian transition functors and then taking the limit in the fiber. 
    
    Note that this limit always exists because all the $\Map_{\Ical}(i, j)$ are finite spaces and since the $F(j)$ are stable. Consequently, $\alpha_i$ is a well-defined functor which is right adjoint to $\ev_i$.

    Let us now prove the second part of the statement. By a dual result to Proposition \ref{FunctorsFromUn}, see also the discussion in Remark \ref{ComparisonWithGHN}, the identity functor $\id:\Gamma\to\Gamma\simeq\mathrm{oplax\;lim}(F)$ corresponds to the section
    $$
        \begin{tikzcd}
            \ev:\Ical^\op\arrow[r] & \UnCart(\Fun(\Gamma, F(-)))
        \end{tikzcd}
    $$
    of the cartesian fibration $\UnCart(\Fun(\Gamma, F(-)))\to\Ical^\op$ mapping $i\in\Ical^\op$ to $\ev_i$. Write $Z_L$ for the full subcategory of sections of $\UnCart(\Fun(\Gamma, F(-)))$ spanned by pointwise left adjoints and $Z_R$ similarly with pointwise right adjoints in $\UnCart(\Fun(F(-), \Gamma))$; then, by taking suitable mapping categories in Theorem B of \cite{HaugsengHebestreitLinskensNuiten}, we deduce that there is an equivalence $Z_L^\op\simeq Z_R$ which is described pointwise by passing to the adjoints. 
    
    In particular, since $\ev$ is in $Z_L$, it corresponds to a point in $Z_R^\op$, i.e. a section
    $$
        \begin{tikzcd}
            \alpha:\Ical\arrow[r] & \UnCart(\Fun(F(-), \Gamma))
        \end{tikzcd}
    $$
    sending $i\in\Ical$ to $\alpha_i$, which gives rises to $\alpha:\Un(F)\to\Gamma$ by Proposition \ref{FunctorsFromUn}.
\end{proof}

When $\Ical$ is further finite, $\alpha$ endows $\Fun_{\catop{\Ical}}(\catop{\Ical}, \UnCart(F))$ with the universal property of $\StabUn(F)$.

\begin{lmm} \label{StronglyFiniteImpliesGeneration}
    Suppose $\Ical$ is strongly finite, then $\Fun_{\catop{\Ical}}(\catop{\Ical}, \UnCart(F))$ is generated under finite limits and retracts (so in particular under small limits) by the images of the $\alpha_i$'s.
\end{lmm}
\begin{proof}
    Let us write $\Gamma$ for $\Fun_{\catop{\Ical}}(\catop{\Ical}, \UnCart(F))$ and $\widehat{\Gamma}:=\Fun_{\Ical^\op}(\Ical^\op, \UnCart(\Pro(F)))$. Remark that the $\ev_i:\widehat{\Gamma}\to\Pro(F(i))$ also have right adjoints $\widehat{\alpha_i}$ given by relative right Kan extension and, using that mapping spaces in $\Ical$ are finite, the explicit formula implies that the following square actually commutes:
    $$
        \begin{tikzcd}
            \Gamma\arrow[r] & \widehat{\Gamma} \\
            F(i)\arrow[r]\arrow[u, "\alpha_i"] & \Pro(F(i))\arrow[u, "\widehat{\alpha_i}"]
        \end{tikzcd}.
    $$

    Note that $\widehat{\Gamma}$ is complete so that, since the $\ev_i$ jointly detect equivalences, it is generated under limits by the images of the $\widehat{\alpha_i}$, see \cite[Corollary 2.5]{Yanovski} and the following remark or the opposite of \cite[Proposition 2.1.2]{HermKIV} for a modern reference. Since $\Pro(F(i))$ is generated under (cofiltered) limits by $F(i)$, it follows that $\widehat{\Gamma}$ is generated under limits by the images of the $\alpha_i$. 
    
    We now claim that $\Gamma$ is contained in the category of cocompact objects of $\widehat{\Gamma}$. Let $\phi\in\Gamma$, we have to show that $\Nat(-, \phi)$ sends cofiltered limits to filtered colimits. Since each $\phi(i)$ is cocompact in $\Pro(F(i))$, each $\Map(-, \phi(i))$ has this property and by \cite[\href{https://kerodon.net/tag/03P7}{Tag 03P7}]{Kerodon}, we have naturally in $\psi\in\Gamma$:
    $$
        \Nat(\psi, \phi)\simeq\lim_{[f:i\to j]\in\TwAr(\Ical)}\Map(\psi(i), \phi(j))
    $$
    In particular, since $\Ical$ is strongly finite, $\TwAr(\Ical)$ is finite and using that filtered colimits commute with finite limits of spaces, $\phi$ is indeed cocompact. From this, we deduce that $\Gamma$ is generated under finite limits and retracts by the images of the $\alpha_i$ since this is more generally the case of all cocompacts objects of $\widehat{\Gamma}$.
\end{proof}

\begin{prop} \label{ExplicitStabUnStronglyFinite}
    Let $\Ical$ be strongly finite and $F:\Ical\to\CatEx$ a functor. Then, precomposition along $\alpha$ induces an equivalence
    $$
        \begin{tikzcd}
            \alpha^*:\FunEx(\Fun_{\catop{\Ical}}(\catop{\Ical}, \UnCart(F)), \Ecal)\arrow[r, "{\simeq}"]& \Fun^{\fbwex}(\Un(F), \Ecal)
        \end{tikzcd}
    $$ 
    so that $\alpha$ induces a canonical equivalence $\StabUn(F)\simeq\Fun_{\catop{\Ical}}(\catop{\Ical}, \UnCart(F))$.
\end{prop}
\begin{proof}
    Once again we write $\Gamma$ for $\Fun_{\catop{\Ical}}(\catop{\Ical}, \UnCart(F))$. First note that $\alpha^*$ has a right adjoint, which is given by the right Kan extension functor $\alpha_*$; to see this, it suffices to argue that if $f:\Un(F)\to\Ecal$ is fiberwise exact, $\alpha_*(f):\Gamma\to\Ecal$ is exact. By applying Lemma \ref{StronglyFiniteImpliesGeneration} to $\FunEx(F(-), \Ecal)$ instead of $F$, we get that $\Fun^{\fbwex}(\Un(F), \Ecal)$ is generated under limits by the images of the right Kan extensions of exact functors $F(i)\to\Ecal$ along $\beta_i:F(i)\to\Un(F)$ so that it suffices to check this for $f$ right Kan extended from some exact $f_i:F(i)\to\Ecal$. But then, $\alpha_*(f)$ is equivalently computed as $(\alpha_i)_*(f_i)$ for $f_i:F(i)\to\Ecal$ and $\alpha_i$ is exact so it must hold.
    
    We now show that the counit map
    $$
        \begin{tikzcd}
            \alpha^*\alpha_*(f)\arrow[r] & f
        \end{tikzcd}
    $$
    is an equivalence. Note that by the same arguments as before, it suffices to prove it for $f:\Un(F)\to\Ecal$ right Kan extended from an exact $f_i:F(i)\to\Ecal$. But then, $\alpha_*(f)\simeq(\alpha_i)_*f_i$ and the right Kan extension along $\alpha_i$ is performed by precomposition by $\ev_i$ by the previous lemma. Therefore, we are reduced to showing that, if $X\in F(j)$
    $$
        \begin{tikzcd}
            f_i(\alpha_j(X)(i))\arrow[r] & f(X)\simeq ((\beta_i)_*f_i)(X)
        \end{tikzcd}
    $$
    where $\beta_i:F(i)\to\Un(F)$ is the inclusion of the fiber. It is a standard fact that right Kan extension along $\beta_i$ is given by the limit
    $$
        \lim_{\gamma\in\Map_{\Ical}(j, i)}f_i(F(\gamma)(X)).
    $$
    The formula for $\alpha_j(X)(i)$ shows that the counit is therefore the limit-comparison map for $f_i$; but each $\Map_{\Ical}(j, i)$ is finite and $f_i$ is exact, so that the map is an equivalence. \\
    
    In particular, this implies that $\alpha^*$ is a localisation (and $\alpha_*$ fully-faithful). We are thus reduced to showing that $\alpha^*$ is conservative. 
    
    By Lemma \ref{StronglyFiniteImpliesGeneration}, we know that $\Gamma$ is generated under finite limits and retracts by the image of $\alpha$. Therefore, if $\eta:f\to g$ is a natural transformation of exact functors $\Gamma\to\Ecal$ which is an equivalence on every $\alpha(X)$, then $\eta$ is also an equivalence on every finite limits of those. Since equivalences are closed under retracts, we deduce that $\eta$ is an equivalence on every object of $\Gamma$, therefore a natural equivalence. Hence, $\alpha^*$ is conservative and this concludes the proof.
\end{proof}

\begin{rmq} \label{2CategoricalHigherSemiAdd}
    By \cite[Theorem 7.1]{GepnerHaugsengNikolaus}, Proposition \ref{ExplicitStabUnStronglyFinite} actually identifies the oplax limit and the oplax colimit of a diagram $F:\Ical\to\CatEx$ whose source is strongly-finite. This is an instance of a 2-categorical ambidexterity phenomenon, more precisely of higher semiadditivity in the sense of \cite[Definition 4.4.2]{HopkinsLurie}.

    A similar phenomenon is investigated in \cite{ChristDycherhoffWaldeLaxAdd}, but in the realm of $\PrExCat{L}$ of (large) presentable stable categories, where the hypothesis of strong-finiteness is therefore no longer necessary.
\end{rmq}

Thanks to this identification, we can now prove the following result, which shows that oplax colimits of non-degenerate hermitian categories indexed by strongly finite diagrams are computed as in $\CatH$.

\begin{prop}\label{UnPNonDegenerateAndFormula}
    Suppose $F:\Ical\to\CatH$ is a functor such that $\Ical$ is strongly finite and every $F(i)$ is a non-degenerate hermitian category. Then the quadratic functor $\Qoppa_F$ is given by
    $$
        \Qoppa_F(\phi)\simeq\colim_{i\in\Ical}\Qoppa_{F(i)}(\phi(i)).
    $$
    Moreover, $\Qoppa_F$ is non-degenerate and the duality $D_{\Qoppa_F}$ is given by
    $$
        D_{\Qoppa_F}(\phi)(j)\coloneqq\colim_{i\in\Ical}\lim_{\gamma\in\Map_{\Ical}(i, j)}F(\gamma)\left[D_{\Qoppa_{F(i)}}(\phi(i))\right].
    $$
    In particular, if every $F(\gamma)$ is duality-preserving, the formula becomes:
    $$
        D_{\Qoppa_F}(\phi)(j)\coloneqq\colim_{i\in\Ical}\lim_{\gamma\in\Map_{\Ical}(i, j)}D_{\Qoppa_{F(j)}}(F(\gamma)(\phi(i))).
    $$
\end{prop}
\begin{proof}
    By Lemma \ref{AlphaiasAdjoints} and Proposition \ref{ExplicitStabUnStronglyFinite}, $\alpha_i:F(i)\to\StabUn(F)$ can be realized as a right adjoint to $\ev_i$. Hence, the left Kan extension functor $(\alpha_i^\op)_!$ is simply given by precomposition by $\ev_i^\op$, which yields the formula for $\Qoppa_F$.
    
    Recall from Remark \ref{IndDualityAsSection} that, as a functor $\StabUn(F)^\op\to\Ind(\StabUn(F))$, the duality is given by
    $$
        D_{\Qoppa_F}(X)\coloneqq\colim_{i\in\Ical}(\alpha_i^\op)_!(\Ind(\alpha_i)\circ D_{\Qoppa_{F(i)}})(X).
    $$
    Under the identification of Proposition \ref{ExplicitStabUnStronglyFinite}, letting $X$ be a section $\phi:\Ical^\op\to\UnCart(F)$ and since $(\alpha_i^\op)_! \simeq (\ev_i^\op)^*$, the formula simplifies to
    $$
        D_{\Qoppa_F}(\phi)\coloneqq\colim_{i\in\Ical}(\Ind(
        \alpha_i)\circ D_{\Qoppa_{F(i)}}\phi(i)).
    $$
    Since $\Ical$ is finite, this colimit can be taken in either $\Ind(\StabUn(F))$ or $\StabUn(F)$ (up to replacing $\Ind(\alpha_i)$ with $\alpha_i$) without altering the result, hence the non-degeneration. Finally, plugging the formula for $\alpha_i$ of Lemma \ref{AlphaiasAdjoints}, we complete the proof.
\end{proof}

\begin{warn} \label{NotPoincaréInGeneral}
    Let $\Ical=\Delta^1$ and $F:\Delta^1\to\CatP$ constant equal to $(\Ccal, \Qoppa)$. Then already $\UnP(F)$ need not be Poincaré. Indeed, it follows from the above formula that the duality sends a section $\phi$ to:
    $$
        D_{\Qoppa_{F}}(\phi)\coloneqq\left[0\longleftarrow D_{\Qoppa}(\phi(1))\right]
    $$
    and this association is not invertible.
\end{warn}

\begin{rmq}
    We note that the higher semiadditivity statement of Remark \ref{2CategoricalHigherSemiAdd} does not extend to $\CatP$ or $\CatH$: the quadratic functors on the oplax limit can be checked to be the limit
    $$
        \lim_{i\in\Ical}\;(\ev_i)^*\Qoppa_{F(i)}
    $$
    which is not in general the same as the colimit of the same functor.
\end{rmq}

The following is also a straightforward consequence from the identification of Proposition \ref{ExplicitStabUnStronglyFinite}.

\begin{prop} \label{IdempotentCompleteStabUn}
    Let $F:\Ical\to\CatEx$ a functor from a strongly-finite $\Ical$ such that each $F(i)$ is idempotent-complete. Then, $\StabUn(F)$ is idempotent-complete as well.
\end{prop}
\begin{proof}
    Using Proposition \ref{ExplicitStabUnStronglyFinite}, we check that $\Fun_{\Ical^\op}(\Ical^\op, \UnCart(F))$ is idempotent-complete. Note that the fully-faithful composite
    $$
        \begin{tikzcd}
            \Fun_{\Ical^\op}(\Ical^\op, \UnCart(F))\arrow[r] & \Fun(\Ical^\op, \UnCart(F)) \arrow[r]& \Fun(\Ical^\op, \Idem(\UnCart(F)))
        \end{tikzcd}
    $$
    preserves and detects splitting idempotents since those are universal limits and these exist and can be computed pointwise at the target. Since each $F(i)$ is idempotent complete, any such splitting $\phi$ is such that $\phi(i)\in F(i)$ and, consequently, is a section of $\UnCart(F)\to\Ical^\op$ as wanted; this concludes.
\end{proof}

Let us now make the following observation: if $f:(\Ccal, \Qoppa)\to (\Dcal, \Qoppa')$ is a duality-preserving functor between non-degenerate hermitian categories, and $(\Ccal, \Qoppa)$ is Poincaré, then for every object $X\in\Ccal$, the map
$$
    \begin{tikzcd}
        f(X)\arrow[r, "{\simeq}"]& D_{\Qoppa'}^\op D_{\Qoppa'}(f(X))
    \end{tikzcd}
$$
is an equivalence, because $\begin{tikzcd}[cramped]X\arrow[r, "{\simeq}"] & D_{\Qoppa}^\op D_{\Qoppa}(X)\end{tikzcd}$ is an equivalence and $f$ is duality-preserving. In particular, $D_{\Qoppa'}$ is fully-faithful on the image of $f$. Now, a non-degenerate hermitian category $(\Dcal, \Qoppa')$ is Poincaré as soon as the duality is fully-faithful since the duality is automatically adjoint to its opposite.

In particular, to show such a $(\Dcal, \Qoppa')$ is Poincaré, it suffices to find duality-preserving functors $f_i:(\Ccal_i, \Qoppa_i)\to (\Dcal, \Qoppa')$ such that the images of the $f_i$'s jointly generate $\Dcal$ under finite colimits. Note that since we found in Remark \ref{NotPoincaréInGeneral} instances where $\UnH(F)$ was not Poincaré, this implies that even if $\Ical$ is strongly finite, $\alpha_i:F(i)\to\UnP(F)$ need not be duality-preserving and in fact this must fail already for $\Ical=\Delta^1$. \\

Let us now investigate more generally when functors of the form $\UnP(G\circ\alpha)\to\UnP(G)$ are duality-preserving. Note that for every $\gamma:\Ical\to\Jcal$ and every $G:\Jcal\to\CatP$, there is a hermitian functor
$$
    \Phi^{\mathrm{p}}_{G, \gamma}:\UnP(G\circ\gamma)\longrightarrow \UnP(G)
$$
induced by the universal properties as oplax colimits of both sides. The first observation is as follows:

\begin{lmm} \label{CofinalInducesCocartesianLift}
    Let $G:\Jcal\to\CatP$ and $\gamma:\Ical\to\Jcal$, and denote $F\coloneqq G\circ\gamma$. Suppose that $\gamma$ is cofinal, then the canonical hermitian functor
    $$
        \Phi^{\mathrm{p}}_{G, \gamma}:\UnP(G \circ \gamma)\longrightarrow \UnP(G)
    $$
    is a cocartesian lift of the functor $\Phi^{\Ex}_{G , \gamma}:\StabUn(F)\longrightarrow\StabUn(G)$, namely the quadratic functor on $\UnP(G)$ is left Kan extended from the one on $\UnP(F)$. 
\end{lmm}
\begin{proof}
    Write $\alpha_i^F:F(i)\to\StabUn(F)$ and $\alpha_j^G:G(j)\to\StabUn(G)$ for the canonical functors. By construction, the composite of $\Phi^{\Ex}_{G, \gamma}\circ\alpha_i^F$ is equivalent to $\alpha_{\gamma(i)}^G$. Consequently, the following square is commutative
    $$
        \begin{tikzcd}
            \Ical\arrow[r, "\Qcal_{F}"]\arrow[d, "\gamma"] & \Fun^q(\StabUn(F))\arrow[d, "(\Phi^{\Ex}_{F, \gamma})_!"] \\
            \Jcal\arrow[r, "\Qcal_{G}"] & \Fun^q(\StabUn(G))
        \end{tikzcd}.
    $$
    where $\Qcal_{F}$ is the functor whose colimit defines $\Qoppa_F$ as in Definition \ref{definition-Qoppa-Ex}. It now follows formally that the natural transformation of quadratic functors given by $\Phi^{\mathrm{p}}_{G, \gamma}$ is the map
    $$
        (\Phi^{\Ex}_{G, \gamma})_!\Qoppa_{G\circ\gamma}\simeq\colim_{i\in\Ical}\Qcal_{G}(\gamma(i))\longrightarrow\colim_{j\in\Jcal}\Qcal_G(j)\simeq\Qoppa_G
    $$
    induced by restricting the diagram of the colimit. If $\gamma$ is cofinal, then this map is an equivalence which concludes.
\end{proof}

Note that in the case where $G$ is constant and $\gamma$ a map of posets, the previous statement is essentially contained in \cite[Proposition 1.5.3(ii)]{HermKII}. Let us also record the following generalization of Lemma \ref{AlphaiasAdjoints}.

\begin{lmm} \label{UpwardsClosedInclusionAdjoint}
    Let $\gamma:\Ical\to\Jcal$ be a fully-faithful functor between strongly-finite categories. Suppose $G \from \Jcal\to\CatEx$ is a functor. Then the functor
    \[
        \Phi^{\Ex}_{G, \gamma}:\StabUn(G\circ\gamma)\longrightarrow\StabUn(G)
    \]
    is fully faithful with a left adjoint.
\end{lmm}
\begin{proof}
    Since both $\Ical$ and $\Jcal$ are strongly finite, it suffices to prove this for the categories of sections in Proposition \ref{ExplicitStabUnStronglyFinite}. In this case, by analogy with Lemma \ref{AlphaiasAdjoints}, we identify $ \Phi^{\Ex}_{G, \gamma}$ as the functor performing the relative right Kan extension in the diagram:
    \[
    \begin{tikzcd}[row sep=large]
    	{\catop{\Ical}}\arrow[r, "\phi"]\arrow[d, "\gamma^\op"] & {\UnCart(G\circ\gamma)}\arrow[r] & {\UnCart(G)}\arrow[d] \\
    	{\catop{\Jcal}}\arrow[rr, equal]\arrow[rru, dotted, "{\Phi^{\Ex}_{G, \gamma}(\phi)}"'] && {\catop{\Jcal}}
    \end{tikzcd}
    \]
    The existence of this relative right Kan extension follows again from \cite[Corollary 4.3.2.14]{HTT} and this functor is automatically right adjoint to precomposition along $\gamma^\op$. The counit $(\Gamma^\op)* \circ  \Phi^{\Ex}_{G, \gamma} \simeq \id$ by the dual of \cite[Proposition 4.3.2.17]{HTT}, and hence we see that $\Phi^{\Ex}_{G, \gamma}$ is fully faithful.
\end{proof}

\begin{defi}
    We say that a functor $\gamma:\Ical\to\Jcal$ is \textit{upwards-closed} if, for every $i\in\Ical$, the induced functor
    $$
        \Ical_{i/}\to\Jcal_{\gamma(i)/}
    $$
    is cofinal.
\end{defi}

If $\Ical$ is a subposet of $\Jcal$ with $\gamma$ the inclusion, then this condition is implied by the more usual definition of upwards-closure, i.e. that if $i\in\Ical$ and $i\leq j$ then $j\in\Ical$.

Our next goal is to generalize Proposition 1.5.3(i) in \cite{HermKII}. For this, let us introduce the following terminology: given a hermitian category $(\Ccal, \Qoppa)$, its $\Ind$-duality is the colimit-preserving functor $D_\Qoppa:\Ind(\Ccal^\op)\to\Ind\Ccal$ canonically associated to $B_\Qoppa$ under the equivalence of categories
$$
    \FunEx(\Ccal^\op\otimes\Ccal^\op, \Sp)\simeq\FunL(\Ind(\Ccal^\op), \Ind(\Ccal))
$$
A hermitian functor $f:(\Ccal, \Qoppa)\to(\Dcal, \Psi)$ is $\Ind$-duality preserving if the usual lax-commuting square involving $\Ind$-dualities and $\Ind(f)$ actually commutes. This recovers precisely the usual notion of duality-preservation in the case where $(\Ccal, \Qoppa)$ and $(\Dcal, \Psi)$ are non-degenerate.

\begin{prop} \label{FibreCofinalInducesCartesianLift}
    Let $\gamma:\Ical\to\Jcal$ be an upwards-closed functor and $G:\Jcal\to\CatH$ a functor. Then, $\Phi^{\mathrm{p}}_{G, \gamma}$ is $\Ind$-duality preserving, i.e. the canonical map 
    $$
        \begin{tikzcd}
            D_{\Qoppa_{G}}\circ\Ind(\Phi^{\Ex}_{G, \gamma})^\op\arrow[r, "\simeq"] & \Ind(\Phi^{\Ex}_{G, \gamma})\circ D_{\Qoppa_{G\circ\gamma}}
        \end{tikzcd}.
    $$ 
    In particular, if $\Ical$ and $\Jcal$ are strongly finite, then by Proposition \ref{UnPNonDegenerateAndFormula}, $\Phi^{\mathrm{p}}_{G, \gamma}$ is duality-preserving. \\
    
    Moreover, supposing further that $\gamma$ is fully-faithful, then the hermitian functor
    $$
        \Phi^{\mathrm{p}}_{G, \gamma}:\UnP(G\circ\gamma)\longrightarrow \UnP(G)
    $$
    is a cartesian lift of the functor $\Phi^{\Ex}_{G, \gamma}:\StabUn(G\circ\gamma)\longrightarrow\StabUn(G)$.
\end{prop}
\begin{proof}
    For a general $F:\Ical\to\CatP$, it is straightforward that:
    $$
        \Qoppa_F\simeq (\alpha_F^\op)_!\left(\colim_{i\in\Ical}(\beta^\op_i)_!\Qoppa_{F(i)}\right)
    $$
    where $\beta_i:F(i)\to\Un(F)$ is the earnest inclusion of the fiber and $\alpha_F:\Un(F)\to\StabUn(F)$ the canonical map. We write $\Qoppa_F^{\mathrm{un}}$ for the functor $\Un(F)^\op\to\Sp$ given by the formula in parenthesis; this functor is such that $\Qoppa_F\simeq(\alpha_F^\op)_!\Qoppa_F^{\mathrm{un}}$. 

    In general, the left Kan extension along $\beta_i^\op$ is computed by the following colimit:
    $$
        (\beta_i^\op)_!\Phi(Y)\simeq \colim_{(F(i)^\op)_{/Y}}\Phi(X)
    $$
    for $\Phi:F(i)^\op\to\Ecal$ any functor whose target admits such colimits and $Y\in F(j)$. Using the cocartesianness of $\Un(F)\to\Ical$ (and the op swapping the direction of the arrow), it follows that this colimit is more simply computed by
    $$
        (\beta_i^\op)_!\Phi(Y)\simeq \colim_{\gamma\in\Map_\Ical(j, i)^\op}\Phi(F(\gamma)(Y))
    $$
    Let $X\in G(\gamma(k))$; the previous discussion shows that $\Qoppa_{G\circ\gamma}^{\mathrm{un}}(X)\to\Qoppa_{G}^{\mathrm{un}}(X)$ is explicitly given by:
    $$
        \colim_{i\in\Ical}\colim_{\sigma\in\Map_{\Ical}(k, i)^\op}\Qoppa_{G(\gamma(i))}(G(\gamma(\sigma))(X))\longrightarrow\colim_{j\in\Jcal}\colim_{\tau\in\Map_{\Jcal}(\gamma(k), j)^\op}\Qoppa_{G(i)}(G(\tau)(X))
    $$
    The global colimit on the left hand side is indexed by $\Ical_{k/}$ whereas the global colimit on the right hand side is indexed by $\Jcal_{\gamma(k)/}$. Since $\gamma$ is upwards-closed, the above map is therefore an equivalence. Hence we have $\Qoppa_{G\circ\gamma}^{\mathrm{un}}\simeq\Qoppa_{G}^{\mathrm{un}}\circ\Phi_{G, \gamma}$ where $\Phi_{G, \gamma}$ is the inclusion of $\Un(G\circ\gamma)\to\Un(G)$.

    Note that a similar argument using the formula of \ref{IndDualityAsSection} also implies that there is an equivalence $B_{\Qoppa_{G\circ\gamma}}^{\mathrm{un}}\simeq B_{\Qoppa_{G}}^{\mathrm{un}}\circ(\Phi_{G, \gamma}^\op\times\Phi_{G, \gamma}^\op)$. \\
    
    We have a commutative square
    $$
        \begin{tikzcd}
            \Un(G\circ\gamma)\arrow[r, "{\Phi_{G, \gamma}}"]\arrow[d, "\alpha_{G\circ\gamma}"] & \Un(G)\arrow[d, "\alpha_{G}"] \\
            \StabUn(G\circ\gamma)\arrow[r, "{\Phi^{\Ex}_{G, \gamma}}"] & \StabUn(G)
        \end{tikzcd}.
    $$
    Hence the transformation $\Qoppa_{G\circ\gamma}\longrightarrow(\Phi^{\Ex}_{G, \gamma})^*\Qoppa_{G}$ is equivalently given by precomposing by $\Qoppa_{G}^{\mathrm{un}}$ the (opposite of) the following Beck-Chevalley transformation:
    \begin{equation*} \tag{B--C} \label{BeckChevalley}
        (\alpha_{G\circ\gamma})_!(\Phi_{G, \gamma})^*\longrightarrow(\Phi^{\Ex}_{G, \gamma})^*(\alpha_G)_!
    \end{equation*}
    associated to the above square when passing to precomposition functors, and this also holds for the $B_\Qoppa$'s in two variables and therefore also for the $\Ind$-dualities. At this point, we already get the first part of the statement: by virtue of Theorem \ref{StabUnDefThm}, the $\alpha_G^*$ are invertible when restricted to exact functors, and squares with invertible vertical legs automatically satisfy Beck-Chevalley conditions. \\
    
    Suppose now $\gamma:\Ical\to\Jcal$ is fully-faithful. Write $\beta_j:G(j)\to\Un(G)$ for the inclusion of the fiber; since all of the functors in the transformation (\ref{BeckChevalley}) commute with colimits, it suffices for our purposes to check that the following transformation is an equivalence
    $$
        (\alpha_{G\circ\gamma})_!(\Phi_{G, \gamma})^*(\beta_j)_!\longrightarrow(\Phi^{\Ex}_{G, \gamma})^*(\alpha_G)_!(\beta_j)_!.
    $$
    Since $\gamma$ is fully faihtful, the functor $\Phi_{G, \gamma}$ is also fully-faithful. It follows from Lemma \ref{UpwardsClosedInclusionAdjoint} that $\Phi^{\Ex}_{G, \gamma}$ is also fully-faithful. But in a square whose horizontal legs are fully-faithful, the Beck-Chevalley condition is always verified for functors left Kan extended from the top left corner. As a colimit of functors left Kan extended from fibers, this is also true for $\Qoppa_{G}^{\mathrm{un}}$, which concludes the proof.
\end{proof}

\subsection{Simplicial complexes and oplax colimits}

\hspace{1.2em} The goal of this subsection is to give a positive example of non-trivial categories $\Ical$ such that for \textit{all} functors $F:\Ical\to\CatP$ the hermitian category $\UnP(F)$ is Poincaré. First, recall the following classical definition:

\begin{defi}
    A \textit{simplicial complex} is a set $K$ such that if $Y\subset X$ is non-empty and $X\in K$ then $Y\in K$. We say that $K$ is a \textit{finite} simplicial complex if every $X\in K$ is a finite set and $K$ itself is finite as a set.

    Given a simplicial complex $K$, we write $\Face(K)$ for the poset of its subsets ordered by inclusion.
\end{defi}

Remark that the space $|\Face(K)|$, obtained by universally inverting all of the arrows of $\Face(K)$, is equivalent to the realization $|K|$ of $K$. Note also that, if $K$ is a finite simplicial complex, then $\Face(K)$ is a strongly finite category. 

\begin{ex}
    Write $[n]=\{0, ..., n\}$ for $n\in\N$. Then, the subsets of $[n]$ form a simplicial complex whose poset of faces is $\Face([n])\simeq\Pcal_{\neq\emptyset}([n])$, i.e. a cube of dimension $n$ with the top vertex removed. This is often called the \textit{barycentric subdivision} of $[n]$. \\
    
    We will write $\CubePunct{n}$ for the \textbf{opposite} of $\Face([n])$. We draw $\CubePunct{2}$ below:
    $$
    \begin{tikzcd}
    	&\{0, 1, 2\}\arrow[rr]\arrow[dl]\arrow[dd] && \{0, 2\}\arrow[dl]\arrow[dd] \\
        \{0, 1\}\arrow[dd]\arrow[rr, crossing over] && \{0\} \\
    	&\{1, 2\}\arrow[rr]\arrow[dl] && \{2\} \\
        \{1\}
    \end{tikzcd}.
    $$
\end{ex}

We will often reduce to finite simplicial complexes, during arguments; our main result to do so is the following Lemma:

\begin{lmm} \label{GeneratedUnderFilteredColimitsByFiniteCplx}
    Let $K$ be a simplicial complex and write $K_\alpha$ for the collection of its finite subcomplexes. Then, $(K_\alpha)$ is a filtered poset under inclusion and we have
    $$
        \begin{tikzcd}\colim\limits_{\alpha}\Face(K_\alpha)^\op\arrow[r, "\simeq"] & \Face(K)^\op\end{tikzcd}
    $$
\end{lmm}
\begin{proof}
    Simplicial complexes are non-empty, so there is at least one finite subcomplex, and finite subcomplexes are stable under union, therefore $(K_\alpha)$ is filtered. Furthermore, any face belongs to a finite simplicial complex, so the inclusions $\Face(K_\alpha)^\op\to\Face(K)^\op$ combine to give
    $$
        \begin{tikzcd}\colim\limits_{\alpha}\Face(K_\alpha)^\op\arrow[r, "\simeq"] & \Face(K)^\op\end{tikzcd}
    $$
    which is both fully-faithful and essentially surjective, hence an equivalence.
\end{proof}

\begin{lmm} \label{OpFacesIsUpwardsClosed}
    Let $K\subset L$ be an inclusion of a subsimplicial complex. Then, the induced functor $\Face(K)^\op\to\Face(L)^\op$ is upwards-closed.
\end{lmm}
\begin{proof}
    Since $\Face(K)^\op$ and $\Face(L)^\op$ are posets, it suffices to see that every face $X\in\Face(L)$ contained in a face $Y\in\Face(K)$ is necessarily a face of $K$; but this is precisely the condition in the definition of a simplicial complex.
\end{proof}

The original idea behind the following lemma dates back to Lurie in \cite[Proposition 3, Lecture 19]{LurieNotesLtheory} and the constant case is contained in \cite[Proposition 6.6.1]{HermKI}.

\begin{lmm} \label{UnPPoincaréCubePunctN}
    Let $F:\CubePunct{n}\to\CatP$ be any functor. Then, $\UnP(F)$ is Poincaré.
\end{lmm}
\begin{proof}
    Since $\CubePunct{n}$ is strongly-finite, we know that the duality on $\UnP(F)$ is non-degenerate. We proceed by induction on $n$; the case $n=0$ being immediate. We remark that every strict subset $T$ of $[n]$ is in the image of a subsimplicial complex $[n-1]\subset[n]$ so in the image of the associated upwards-closed inclusion $\beta:\CubePunct{n-1}\to\CubePunct{n}$.

    By Proposition \ref{FibreCofinalInducesCartesianLift}, the induced functor $\UnP(F\circ\beta)\to\UnP(F)$ is duality-preserving and the induction hypothesis guarantees that its source is a Poincaré category. To conclude, it thus suffices to see that for $X\in F([n])$, the map
    $$
        \begin{tikzcd}
            \alpha_{[n]}(X)\arrow[r, "{\simeq}"] & D_{\Qoppa_F}D_{\Qoppa_F}^\op\alpha_{[n]}(X)
        \end{tikzcd}
    $$
    is an equivalence.
    
    Using the explicit formula for the duality of Proposition \ref{UnPNonDegenerateAndFormula}, we will prove the following:
    \begin{align*}
        D_{\Qoppa_F}(\alpha_{[n]}(X))([n])&\simeq \Sigma^n D_{\Qoppa_{F([n])}}(X) \\
        \text{For }S\neq[n],\text{ }D_{\Qoppa_F}(\alpha_{[n]}(X))(S)&\simeq 0.
    \end{align*}
    In both cases, the formula of Proposition \ref{UnPNonDegenerateAndFormula} is the colimit of a punctured cube. To compute this colimit, we will invoke usual arguments about total cofibers of cubes: namely, recall that given a non-necessarily cocartesian cube $X$, its total cofiber is the cofiber map
    $$
        \totcofib(X) :=\cofib\left(\colim_{S\in\CubePunct{n}}X(S)\to X(\emptyset)\right).
    $$
    In our cases, we can always complete the cube by setting $X(\emptyset)=0$ so that computing the total cofiber is simply computing a shift of the colimit we are interested in. A standard fact about total cofibers (see \cite{Goodwillie2}) is that they can be computed iteratively: namely, choosing a direction to view a $n$-cube as a map of $(n-1)$-cube, the total cofiber is also the cofiber of the map between the total cofibers of the restricted $(n-1)$-cubes:
    $$
        \totcofib(X)\simeq\cofib\left(\totcofib(X_{\textup{left}})\to\totcofib(X_{\textup{right}})\right).
    $$
    In particular, we deduce from this that a map between cocartesian $(n-1)$-cubes is a cocartesian $n$-cube. Note that, in the second case, the formula asks to compute the colimit of a punctured cube where the subdiagram of those subsets $T$ such that $T\supseteq S$ --- which spans a cube of size $|S|$ --- is constant equal to $D_{\Qoppa_{F(S)}}(F(\gamma)(X))$ where $\gamma:[n]\to S$, with transitions maps being the identity and the rest of the cube is zero. In particular, thanks to the preceding discussion about total cofibers, we see that the diagram with an added zero at the last vertex is cocartesian, hence the second equivalence.
    
    For the first case, we are trying to compute the colimit of a punctured cube whose top vertex is $D_{\Qoppa_{F([n])}}(X)$ and every other vertex is zero (by virtue of having no maps $S\to [n]$ since $[n]$ is not included into any strict subset). A straightforward induction using total cofibers yields the fact that this is $\Sigma^n D_{\Qoppa_{F([n])}}(X)$ as wanted. \\

    To conclude, it suffices to remark that for every non-empty $S\subset[n]$, the colimit formula of Proposition \ref{UnPNonDegenerateAndFormula} applied to compute $D_{\Qoppa_F}(D_{\Qoppa_F}^\op\alpha_{[n]}(X))(S)$ is always of the first case, i.e. every vertex is zero except the top one, which is given by 
    $$
        D_{\Qoppa_{F([n])}}(\Sigma^n F(\gamma)(D_{\Qoppa_{F([n])}}(X)))
    $$ 
    where $\gamma:[n]\to S$. Since $D_{\Qoppa_{F([n])}}$ is contravariant, it sends $\Sigma^n$ to $\Omega^n$, which cancels out with the suspension appearing after taking the colimit of the cube. Consequently, we get that
    $$
        D_{\Qoppa_F}(D_{\Qoppa_F}^\op\alpha_{[n]}(X))(S)\simeq F(\gamma)(X)
    $$
    the right hand side being precisely $\alpha_{[n]}(X)(S)$ by Lemma \ref{AlphaiasAdjoints}. 

    Let us finally check that we have indeed proven that the correct map is an equivalence. Unraveling the formulas, such map is precisely the canonical map:
    $$
        \begin{tikzcd}
            \alpha_{[n]}(X)(j)\arrow[r] & \colim\limits_{i\in\CubePunct{n}}\lim\limits_{\gamma\in\Map(i, j)}\lim\limits_{k\in\CubePunct{n}}\colim\limits_{\delta\in\Map(k, i)}D_{\Qoppa_{F(j)}} D_{\Qoppa_{F(j)}}^\op F(\gamma)F(\delta)(\alpha_{[n]}(X)(k)).
        \end{tikzcd}
    $$
    Using the natural equivalence $D_{\Qoppa_{F(j)}} D_{\Qoppa_{F(j)}}^\op\longrightarrow\id$, this amounts to showing that the map
    $$
        \begin{tikzcd}
            \alpha_{[n]}(X)(j)\arrow[r] & \colim\limits_{i\in\CubePunct{n}}\lim\limits_{\gamma\in\Map(i, j)}\lim\limits_{k\in\CubePunct{n}}\colim\limits_{\delta\in\Map(k, i)}F(\gamma)F(\delta)(\alpha_{[n]}(X)(k))
        \end{tikzcd}
    $$
    is an equivalence. Note that the mapping spaces appearing in the above formula are either contractible (if say $i\leq j$ or $k\leq i$) or empty, hence we can write the map once again as
    $$
        \begin{tikzcd}
            \alpha_{[n]}(X)(j)\arrow[r] & \colim\limits_{i\in\CubePunct{n}}\lim\limits_{k\in\CubePunct{n}}\begin{cases*}
                F(j\leq k)(\alpha_{[n]}(X)(k)) \text{ if }j\leq i \leq k \\ 0 \text{ otherwise}
            \end{cases*}
        \end{tikzcd}.
    $$
    Up to commuting the finite limit and the finite colimit, the above argument is therefore precisely computing the right hand side by checking that the canonical map from $\alpha_{[n]}(X)(j)$ is an equivalence. This concludes.
\end{proof}

\begin{rmq} \label{AlphaNotPoincareButEvYes}
    In particular, the formulas we produced during the proof show that even though $\UnP(F)$ is Poincaré, the functor $\alpha_{[n]}:F([n])\to\UnP(F)$ is \textit{not} a Poincaré functor (unless $n=0$).

    However, the computation of the duality at the point $[n]$ extends to general sections $\phi$, so we have shown that $\ev_{([n])}$ can be refined to a Poincaré functor
    $$
        \begin{tikzcd}
            \ev_{([n])}:\UnP(F)\arrow[r] & (F([n]), \Sigma^n\Qoppa_{F([n])})
        \end{tikzcd}
    $$
    provided we shift the quadratic functor on either side.
\end{rmq}

\begin{cor} \label{UnPPoincareFiniteSimplicialComplex}
    Let $K$ be a  simplicial complex and denote $I$ the opposite of its poset of faces. Then, for every functor $F:I\to\CatP$, $\UnP(F)$ is Poincaré.
\end{cor}
\begin{proof}
    By definition, a simplicial complex can be exhausted by maps from the simplicial complex $\Delta^n$. Those maps induce upwards-closed inclusions when passing to the opposite of the poset of faces so that Proposition \ref{FibreCofinalInducesCartesianLift} implies that the functor
    $$
        \begin{tikzcd}
            \UnP(F\mid_{\CubePunct{n}})\arrow[r] & \UnP(F)
        \end{tikzcd}
    $$
    is duality-preserving; since every point in $\StabUn(F)$ is in the image of one such map, the result follows from Lemma \ref{UnPPoincaréCubePunctN}.
\end{proof}

\section{Assembly and simplicial complexes}
\subsection{Localizations and oplax colimits}

\hspace{1.2em} For this section, we fix a functor $i:\Ical\to\Jcal$ between categories and $F:\Jcal\to\CatP$ a functor. We are interested in the relation between the oplax colimit of $F$ and $F\circ i$. As we have already remarked, the universal property of oplax colimits induces a hermitian functor:
$$
    \begin{tikzcd}
        \Phi^{\mathrm{p}}_{F, i}:\UnP(F\circ i)\arrow[r] & \UnP(F).
    \end{tikzcd}
$$
By building on work by Hinich \cite{Hinich}, we can say more when $\Ical\to\Jcal$ is a localisation:

\begin{lmm} \label{HinichTrick}
    Suppose $i:\Ical\to\Jcal$ is a localization at $\Wcal$, then both maps
    $$
        \Phi_{F, i}:\Un(F\circ i)\to\Un(F)\text{ and }\Phi_{F, i}^{\Ex}:\StabUn(F\circ i)\to\StabUn(F)
    $$
    are localizations, the former at the collection of arrows $\Wcal'$ of cocartesian lifts of $\Wcal$ and the latter at its image $\overline{\Wcal}$ under $\Un(F\circ i)\to\StabUn(F\circ i)$. 
\end{lmm}
\begin{proof}
    The claim for the first functor is a consequence of \cite[Proposition 2.1.4]{Hinich}. Since we will need to prove the second claim and the strategy will be the same, let us present a different proof of this first fact.
    
    Let $G:\Jcal\to\CatInfty$ be any functor, then the (equivalent) horizontal arrows:
    $$
        \begin{tikzcd}
            \Fun_{\Jcal}(\Jcal, \UnCart(G))\arrow[d, "="]\arrow[r]& \Fun_{\Ical}(\Ical, \UnCart(G\circ i))\arrow[d, "\simeq"]\\
            \Fun_{/\Jcal}(\Jcal, \UnCart(G))\arrow[r, "i^*"]& \Fun_{/\Jcal}(\Ical, \UnCart(G))
        \end{tikzcd}
    $$
    are fully-faithful with image those functors $\Ical\to\UnCart(F)$ over $\Jcal$ which factor through $\Jcal$, i.e. which invert $\Wcal$. Leveraging Proposition \ref{FunctorsFromUn} and applying the above to $G\coloneqq\Fun(F(-), \Ecal)$ for some $\Ecal$, we deduce that the functor
    $$
        (\Phi_{F, i})^*:\Fun(\Un(F), \Ecal)\longrightarrow\Fun(\Un(F\circ i), \Ecal)
    $$
    is fully-faithful with image those functors which invert $\Wcal'$. This concludes for the first claim; we note in addition, that since the functor is given naturally in $\Ecal$, this actually supplies the functor $\Phi_{F, i}$ by fully-faithfulness of the Yoneda embedding. \\

    The second claim follows \textit{mutatis mutandis}, using the first observation for $G\coloneqq\FunEx(F(-), \Ecal))$ instead for a stable $\Ecal$ --- this is enough by \cite[Theorem I.3.3]{NikolausScholze} --- and remarking as previously that this defines the exact functor $\Phi_{F, i}^{\Ex}$. 
\end{proof}

Let us add the following description of the maps inverted by $\Phi^{\Ex}_{F, i}$, which is a folklore fact of stable categories for which we did not find a proof in the literature.

\begin{lmm} \label{GeneratorsKernelOfVerdierProj}
    With the same notations as the previous lemma, the kernel of $\Phi_{F, i}^{\Ex}$ is generated as a full stable subcategory of $\StabUn(F)$ closed under retracts by the fibers of maps in $\overline{\Wcal}$.
\end{lmm}
\begin{proof}
    We first note that it is clear that the kernel of $\Phi_{F, i}^{\Ex}$ contains the fibers of maps in $\overline{\Wcal}$, any of their finite colimits and any retracts of them. We are interested in the converse direction; suppose $\alpha:Z\to Z'$ is a map in $\ker(\Phi_{F, i}^{\Ex})$ such that
    $$
        \begin{tikzcd}
            \alpha^*:\map(Z', \ker(f))\arrow[r, "{\simeq}"] & \map(Z, \ker(f))
        \end{tikzcd}
    $$
    for every $f:X\to Y\in\overline{\Wcal}$. This implies that
    $$
        \begin{tikzcd}
            f_*:\map(\coker(\alpha), X)\arrow[r, "{\simeq}"] & \map(\coker(\alpha), Y)
        \end{tikzcd}
    $$
    which in turns implies that $\coker(\alpha)$ is both local with respect to $\Wcal$ and sent to zero by the localisation (as this is the case for $Z$ and $Z'$). Hence, $\coker(\alpha)\simeq 0$ and, finally, $\alpha:Z\to Z'$ is an equivalence, which shows that the collection $\{\ker(f)\mid f\in\overline{\Wcal}\}$ detects equivalences in the kernel of $\Phi_{F, i}^{\Ex}$. This suffices by A.1.8 of \cite{HermKII}. 
\end{proof}

We can now show the existence of a hermitian functor.

\begin{lmm} \label{CocartesianTransition}
    If $i:\Ical\to\Jcal$ is a localisation and $F:\Jcal\to\CatH$ a functor, then the exact functor
    $\Phi_{F, i}^{\Ex}$ of the previous lemma refines to a functor between hermitian categories 
    $$
        \Phi_{F, i}^{\mathrm{p}}:\UnP(F\circ i)\longrightarrow\UnP(F)
    $$
    which is a cocartesian lift of $\Phi_{F, i}^{\Ex}$ with respect to the cocartesian fibration $\CatH\to\CatEx$.
\end{lmm}
\begin{proof}
    We proceed as in the previous lemma. Let $G\coloneqq\FunH(F(-), (\Ccal, \Qoppa))$ for $(\Ccal, \Qoppa)$ a Poincaré category, then
    $$
        \begin{tikzcd}
            \Fun_{/\Jcal}(\Jcal, \UnCart(G))\arrow[r, "i^*"]& \Fun_{/\Jcal}(\Ical, \UnCart(G))
        \end{tikzcd}
    $$
    is fully-faithful with image functors inverting $\Wcal$. Hence, there is a fully-faithful functor
    $$
        \begin{tikzcd}
            \FunH(\UnP(F), (\Ccal, \Qoppa))\arrow[r] & \FunH(\UnP(F\circ i), (\Ccal, \Qoppa))
        \end{tikzcd}
    $$
    which is natural in $(\Ccal, \Qoppa)$. We deduce that there is a hermitian functor
    $$
        \begin{tikzcd}
            \Phi_{F, i}^{\mathrm{p}}:\UnP(F\circ i)\arrow[r] & \UnP(F)
        \end{tikzcd}.
    $$
    Since the following diagram commutes, $\Phi_{F, i}^{\mathrm{p}}$ recovers $\Phi_{F, i}^{\Ex}$ on underlying categories:
    $$
        \begin{tikzcd}
            \FunH(\UnP(F), (\Ccal, \Qoppa))\arrow[r, "{(\Phi_{F, i}^{\mathrm{p}})^*}"]\arrow[d] & \FunH(\UnP(F\circ i), (\Ccal, \Qoppa))\arrow[d] \\
            \FunEx(\StabUn(F), \Ccal)\arrow[r, "{(\Phi_{F, i}^{\Ex})^*}"] & \FunEx(\StabUn(F\circ i), \Ccal)
        \end{tikzcd}.
    $$
    The final claim follows from checking that this square is actually cartesian; since the vertical maps are left fibrations, it suffice to check that vertical fibers are equivalent. Again using universal properties, this is more explicitly the square:
    $$
        \begin{tikzcd}
            \Fun_{/\Jcal}(\Jcal, \UnCart(\FunH(F(-), (\Ccal, \Qoppa))))\arrow[r]\arrow[d] & \Fun_{/\Jcal}(\Ical, \UnCart(\FunH(F(i(-)), (\Ccal, \Qoppa))))\arrow[d] \\
            \Fun_{/\Jcal}(\Jcal, \UnCart(\FunEx(F(-), \Ccal)))\arrow[r] & \Fun_{/\Jcal}(\Ical, \UnCart(\FunEx(F(i(-)), \Ccal)))
        \end{tikzcd}.
    $$
    On vertical fibers, fixing some section $s:\Jcal\to\UnCart(\FunEx(F(-), \Ccal))$, the induced map is given by
    $$
        \begin{tikzcd}
            \Fun_{/\Jcal}(\Jcal, \UnCart(\Nat((s_{(-)}^\op)_!\Qoppa_{F(-)}, \Qoppa')))\arrow[r] & \Fun_{/\Jcal}(\Ical, \UnCart(\Nat((s_{(-)}^\op)_!\Qoppa_{F(-)}, \Qoppa')))
        \end{tikzcd}.
    $$
    This is also a fully-faithful functor since $i$ is a localisation. Thus it suffices to argue that any functor over $\Jcal$
    $$
        \begin{tikzcd}
            t:\Ical\arrow[r] & \UnCart(\Nat((s_{(-)}^\op)_!\Qoppa_{F(-)}, \Qoppa'))
        \end{tikzcd}
    $$
    necessarily factors through $i$, i.e. inverts $\Wcal$. For any $w \in \Wcal$, the image of $w$ fits in a triangle where in the one arrow consist of a map in the fibre (as is hence automatically an equivalence since we have unstraightened a space valued functor) and the other arrow is a Cartesian lift of the map $i(w)$ and is thus also an equivalence. Hence the image of $w$ is too.
\end{proof}

\begin{rmq}
    By \cite[Proposition 7.1.10]{CisinskiBook}, localisations are cofinal so one could also argue using Lemma \ref{CofinalInducesCocartesianLift}. 
\end{rmq}

In particular, the characterization as a cocartesian lift identifies $\Phi_{F, i}^{\mathrm{p}}$ as the cofiber functor associated to the cartesian lift of $\ker(\Phi_{F, i}^{\Ex})\to\StabUn(F\circ i)$ with target $\UnP(F\circ i)$. The following proposition becomes now within reach:

\begin{prop} \label{PVprojectionLocalization}
    Let $i:\Ical\to\Jcal$ be a localization and $F:\Jcal\to\CatP$ a functor such that $\UnP(F\circ i)$ is Poincaré and the kernel of $\Phi_{F, i}^{\Ex}$ is stable under its duality. Then,
    $$
        \begin{tikzcd}
            \Phi_{F, i}^{\mathrm{p}}:\UnP(F\circ i)\arrow[r] & \UnP(F)
        \end{tikzcd}
    $$
    is a Poincaré--Verdier projection. In particular, $\UnP(F)$ is a Poincaré category.
\end{prop}
\begin{proof}
    The criterion of \cite[Corollary 1.1.6]{HermKII} guarantees that Poincaré--Verdier projections are precisely Poincaré functors which are cocartesian lifts of Verdier projections. Hence the previous two lemmata conclude provided we can show that $\Phi_{F, i}^{\mathrm{p}}$ is a Poincaré functor between Poincaré categories.
    
    The observation preceding the statement of Proposition \ref{PVprojectionLocalization} combined with the fact that $\CatP\to\CatH$ preserves and detects colimits \cite[Proposition 6.1.4]{HermKI} (but see also \cite[Example 1.1.7]{HermKII}) ensures that this is automatic as soon as the cartesian lift of 
    $$
        \begin{tikzcd}
            \ker(\Phi_{F, i}^{\Ex})\arrow[r] & \StabUn(F\circ i)
        \end{tikzcd}
    $$ 
    with target $\UnP(F\circ i)$ is itself a Poincaré functor between Poincaré categories. The category $\UnP(F\circ i)$ is Poincaré by hypothesis. Since the kernel of $\Phi_{F, i}^{\Ex}$ is stable under the duality of $\UnP(F\circ i)$, it is also a Poincaré category when equipped with the restricted quadratic functor and the cartesian lift is automatically a Poincaré functor. This concludes.
\end{proof}

\subsection{The assembly map of bordism-invariant localizing functors}

\hspace{1.2em} Given a simplicial complex $K$, for any $F:|K|\to\CatP$, we have built a hermitian functor
$$
    \begin{tikzcd}
        \Phi^{\mathrm{P}}_{F}\colon\UnP(F|_{\Face(K)^\op})\arrow[r] & \UnP(F)\simeq\colim_{|K|}F
    \end{tikzcd}
$$
which is the cocartesian lift of a localization by Lemma \ref{CocartesianTransition}. Here, we note that we chose to work with the opposite of $\Face(K)$ but by a game of opposites, there is an equivalence
$$
    \UnP(F|_{\Face(K)^\op})\simeq\Fun_{\Face(K)}(\Face(K), \UnCart(F|_{\Face(K)^\op})).
$$
We also know that $\UnP(F|_{\Face(K)^\op})$ is a Poincaré category by Corollary \ref{UnPPoincareFiniteSimplicialComplex}. Since $\CatP$ is closed under colimits, we know that $\colim_{|K|}F$ is also Poincaré.

However, we still do not know whether $\Phi^{\mathrm{P}}_{F}$ is a Poincaré functor. By Proposition \ref{PVprojectionLocalization}, it suffices to prove that its kernel is stable under the duality. We begin by treating the case where $K$ is finite; then, $\Face(K)^\op$ is a strongly-finite category and therefore we have access to a simpler model for $\UnP(F|_{\Face(K)^\op})$ and consequently for the kernel as well.

\begin{prop} \label{KernelStableUnderDualityFiniteCpx}
    Let $K$ be a finite simplicial complex and $F:|K|\to\CatP$ a functor. Then, the kernel $\ker(\Phi^{\mathrm{P}}_{F})$ is stable under the duality of $\UnP(F|_{\Face(K)^\op})$.
\end{prop}
\begin{proof}
    If $K\subset L$ is an inclusion of finite simplicial complexes, it induces an upwards-closed inclusion of their opposite of posets of faces by Lemma \ref{OpFacesIsUpwardsClosed}. In particular, by Proposition \ref{FibreCofinalInducesCartesianLift}, we get for every $F:|L|\to\CatP$ a duality-preserving
    $$
        \begin{tikzcd}
            \UnP(F|_{\Face(K)^\op})\arrow[r] & \UnP(F|_{\Face(L)^\op})
        \end{tikzcd}
    $$
    between Poincaré categories (in fact a split-Verdier inclusion). Moreover, this functor comes with a homotopy making the following square commute
    $$
        \begin{tikzcd}
            \UnP(F|_{\Face(K)^\op})\arrow[r]\arrow[d] & \UnP(F|_{\Face(L)^\op})\arrow[d] \\
            \UnP(F|_{|K|})\arrow[r] & \UnP(F|_{|L|})
        \end{tikzcd}
    $$
    where the bottom row is simply the Poincaré functor induced by taking colimits. In particular, there is an induced functor of hermitian categories
    $$
        \begin{tikzcd}
            \ker(\Phi^\mathrm{P}_{F, K})\arrow[r] & \ker(\Phi^\mathrm{P}_{F, L})
        \end{tikzcd}.
    $$
    We now see that, if the left hand side is stable under the duality of $\UnP(F|_{\Face(K)^\op})$, then its image in the right hand side must be stable under the duality of $\UnP(F|_{\Face(L)^\op})$. Since every simplicial complex can be exhausted by subcomplexes of the form $\Delta^n$, we are reduced to prove the case where $\Face(K)^\op$ is $\CubePunct{n}$, the opposite of the poset of faces of the simplicial complex $[n]$. \\

    We argue by induction on $n$. The case $n=0$ is trivial. By Proposition \ref{GeneratorsKernelOfVerdierProj}, the kernel is generated by the fibers
    $$
        z^{\gamma}_x:=\fib\left(\begin{tikzcd}\alpha_i(x)\arrow[r] & \alpha_j(F(\gamma)(x))\end{tikzcd}\right)
    $$
    for $x\in F(i)$ and $\gamma:i\to j$ a map in $\CubePunct{n}$ and we have to check that $D_\Qoppa(z^{\gamma}_x)$ is still in this kernel. In fact, using the previous discussion, we can further reduce the checks to be made: it suffices to consider the case where $\gamma:i\to j$ is such that $i=[n]$ is the initial point of the punctured cube since any other map can be dealt with through the inclusion of a face $\Delta^{n-1}\subset\Delta^n$ (i.e. by applying the inductive hypothesis to those). Moreover, note that, if $i\xrightarrow{\gamma} j\xrightarrow{\delta} k$ is a composite in $\CubePunct{n}$, there is a short exact sequence
    $$
        \begin{tikzcd}
            z^\gamma_x\arrow[r] & z^{\delta\circ\gamma}_x\arrow[r] & z^\delta_{F(\gamma)x}
        \end{tikzcd}
    $$
    so that it suffices to show that $D_\Qoppa$ sends the extremal terms to the kernel to conclude. The statement holds for the rightmost one by the inductive assumption. Hence, we reduce to treat the case $z^{[n]\to[n-1]}_x$ for the $n+1$ inclusions $\gamma:[n-1]\subset[n]$; by symmetry, we will only explain the case where $\gamma$ includes $[n-1]$ as $\{0, ..., n-1\}$ in $[n]$.
    
    Note that in Lemma \ref{UnPPoincaréCubePunctN} we have already computed the duality of $\alpha_{[n]}(x)$ and found it was zero everywhere except at the point $[n]$, where it is given by $\Sigma^n D_{\Qoppa_{F([n])}}(x)$. This also holds for $[n-1]$ with one less shift since the inclusion $\Delta^{n-1}\to\Delta^n$ induces a duality-preserving, fully-faithful functor on oplax colimits. In particular, since $D_\Qoppa$ is contravariantly exact, we deduce from the fiber defining $z^\gamma_x$ the following:
    $$
         D_{\Qoppa}(z^{\gamma}_x)(S)\simeq\left\{\begin{aligned}
            \Sigma^n D_{\Qoppa_{F([n])}}(x)\text{ if }S=[n] \\
            F([n]\to [n-1])(\Sigma^n D_{\Qoppa_{F([n])}}(x))\text{ if }S=\{0, ..., n-1\} \\
            0\text{ otherwise}
        \end{aligned}\right.
    $$
    with the canonical maps and homotopies between them. Now we remark that, since $\Delta^n$ is contractible and $F$ is restricted along the localization $\CubePunct{n}\to\Delta^n\simeq *$, all the maps $F([n]\leq [n-1])$ are canonically homotopic to the identity. Moreover, the functor $\UnP(F\mid_{\CubePunct{n}})\to F(*)$ sends a section to the right Kan extension along $\CubePunct{n}\to *$, i.e. the limit of the diagram. 
    
    Hence, the above diagram, once completed with a zero at the top vertex, defines a cube whose total fiber is zero, since these can be computed iteratively. In particular, it follows that the limit $\lim_{S\neq\emptyset} D_{\Qoppa}(z^{\gamma}_x)(S)$ vanishes, so $D_{\Qoppa}(z^{\gamma}_x)$ belongs to the kernel as wanted.
\end{proof}

\begin{rmq} \label{FormulaForZGammaX}
    It will be useful in the rest of the section to have a formula for the more general $z^{\delta}_y$. Using Lemma \ref{AlphaiasAdjoints}, we get an explicit formula for the $\alpha_{[k]}$ which is such that if $\delta:[k]\leq [m]$ and $y\in F([k])$:
    $$
         z^{\delta}_y(S)\simeq\left\{\begin{aligned}
            F([k]\leq S)(y)\text{ if }[k]\leq S\text{ but }[m]\not\leq S \\
            0\text{ otherwise}.
        \end{aligned}\right. 
    $$
    We also record the fact used in the proof that, for a general $\gamma:S\to T$, $D_{\Qoppa}(z^{\gamma}_x)$ is supported on all the $U$ such that $S\to U \to T$.
\end{rmq}

By combining Proposition \ref{PVprojectionLocalization} and Proposition \ref{KernelStableUnderDualityFiniteCpx}, we get that:
$$
    \begin{tikzcd}
        \ker(\Phi^{\mathrm{P}}_{F})\arrow[r] & \UnP(F|_{\Face(K)^\op})\arrow[r, "{\Phi^{\mathrm{P}}_{F}}"] & \UnP(F)\simeq\colim_{|\Face(K)^\op|}
    \end{tikzcd}
$$
is a Poincaré-Verdier sequence when $K$ is a finite simplicial complex. But note that Poincaré-Verdier sequences are stable under filtered colimits since finite limits commute with filtered colimits in $\CatP$, hence we can combine Lemma \ref{GeneratedUnderFilteredColimitsByFiniteCplx} and the functoriality of oplax colimits of Lemma \ref{OplaxColimitsPreserveColimits} to get the following result in full generality:

\begin{cor} \label{PVSequenceForOplaxColim}
    Let $K$ be a simplicial complex and $F:|K|\to\CatP$ a functor. Then, the following is a Poincaré-Verdier sequence
    $$
    \begin{tikzcd}
        \ker(\Phi^{\mathrm{P}}_{F})\arrow[r] & \UnP(F|_{\Face(K)^\op})\arrow[r, "{\Phi^{\mathrm{P}}_{F}}"] & \colim_{|K|}F.
    \end{tikzcd}
    $$
\end{cor}

We now claim that, under a bordism-invariant functor and taking $K=\CubePunct{n}$, the left hand side will vanish. To show this, let us first recall that a subcategory $\Lcal\subset\Ccal$ of a Poincaré category $(\Ccal, \Qoppa)$ is isotropic if $\Qoppa$ vanishes on $\Lcal$ and if the inclusion $\Lcal\to\Ccal$ admits a right adjoint. In particular, such an isotropic category induces an inclusion
$$
    \Lcal\subset\Lcal^{\perp}
$$
where $\Lcal^\perp$ is the full subcategory of $\Ccal$ of those $X$ such that $B_\Qoppa(\Lcal, X)\simeq 0$. We say that $\Lcal$ is a Lagrangian if this map is an equivalence. A Poincaré category admitting a Lagrangian category is called metabolic.

\begin{lmm} \label{EasierCriterionMetabolic}
    Suppose $\Lcal$ is a subcategory of $(\Ccal, \Qoppa)$ such that $\Qoppa(\Lcal)=0$ and $\Lcal, D_\Qoppa(\Lcal)$ jointly generate all of $\Ccal$ as a thick subcategory. Then, $\Lcal$ is a Lagrangian for $(\Ccal, \Qoppa)$.
\end{lmm}
\begin{proof}
    The two conditions we have to check are that $\Lcal\subset\Ccal$ has a right adjoint and that $D_\Qoppa:\Lcal\to\Lcal^{\perp}$ is an equivalence. Note first that, for any $X, Y\in\Lcal$, the mapping spectrum
    $$
        \map_\Ccal(X, D_\Qoppa(Y))\simeq B_\Qoppa(X, Y)
    $$
    vanishes since $B_\Qoppa(X, Y)$ is a retract of $\Qoppa(X\oplus Y)\simeq 0$. The two conditions now follow from checking that there exists, for every $X\in\Ccal$, an exact sequence
    $$
        X_l\longrightarrow X\longrightarrow D_\Qoppa(X_d)
    $$
    with $X_l,X_d\in\Lcal$. Indeed, together with the previous vanishing of the mapping spectra, this implies that there is a semi-orthogonal decomposition of $\Ccal$ by $\Lcal$ and $D_\Qoppa(\Lcal)$ so that we have the wanted right adjoint (see $(v)\implies (ii)$ in Proposition 6.2 of \cite{SaunierWinges} for instance). Furthermore, we see that if $X\in\Lcal^{\perp}$ so is $D_\Qoppa(X_d)$ by closure under fibers of $\Lcal^\perp$ and the by fact that $D_\Qoppa$ carries $\Lcal$ to $\Lcal^\perp$. But then $D_\Qoppa(\Lcal)\cap\Lcal^\perp=0$ and thus $X\simeq X_l$ which shows that $X\in\Lcal$.

    We now focus on producing the exact sequence. Let $\Fcal$ be the full subcategory of $\Ccal$ for which there exists such a sequence. Of course, $\Fcal$ contains both $\Lcal$ and $D_\Qoppa(\Lcal)$, thus it suffices to see that $\Fcal$ is stable under direct sums and fibers; the former is straightforward by simply summing the given sequences. Suppose now given $f:X\to Y$ and exact sequences $X_l\to X\to D_\Qoppa(X_d)$ as well as $Y_l\to Y\to D_\Qoppa(Y_l)$. Then, there exist dotted lifts making the following diagram commute:
    $$
        \begin{tikzcd}
            X_l\arrow[r, dotted]\arrow[d] & Y_l\arrow[d] \\
            X\arrow[r, "f"]\arrow[d] & Y\arrow[d] \\
            D_\Qoppa(X_d)\arrow[r, dotted] & D_\Qoppa(Y_d)
        \end{tikzcd}.
    $$
    Indeed, since taking cofibers is functorial, it suffices to see why the top vertical dashed exists; but we already argued that $\Map(\Lcal, D_\Qoppa(\Lcal))=0$. In particular, the composite $X_l\to Y\to D_\Qoppa(Y_d)$ is null, so that it factors through the fiber of the second map, which is $Y_l$.

    Now, the wanted exact sequence for the fiber of $f$ can be obtained by simply taking horizontal fibers in the above diagram. This concludes.
\end{proof}

\begin{lmm} \label{KernelMetabolicCubes}
    Let $(\Ccal, \Qoppa)$ be a idempotent-complete Poincaré category and denote $F:\CubePunct{n}\to\CatP$ the constant functor equal to $(\Ccal, \Qoppa)$. Then, the kernel of the Poincaré functor
    $$
        \UnP(F)\longrightarrow (\Ccal, \Qoppa)
    $$
    is a metabolic Poincaré category.
\end{lmm}
\begin{proof}
    Note that this situation fits in the previous set-up, since $\CubePunct{n}$ is the opposite of the poset of faces of the contractible simplicial complex $[n]$, hence any functor $F:|\Delta^n|\to\CatP$ is just the datum of a Poincaré category.

    Pick a way to view $\CubePunct{n}$ as a map of $(n-1)$-cubes --- for simplicity, let us choose the one that removes the vertex $\{n\}$ so that the maps between the two $(n-1)$-cubes are given pointwise by the edges $\gamma_S:S\cup\{n\}\to S$ for some $S\in\CubePunct{n-1}$. We define a thick subcategory of the kernel as follows:
    $$
        \Lcal:=\langle z^{\gamma_S}_x \mid \gamma_S:S\cup\{n\}\to S\text{ for }S\in\CubePunct{n-1}\rangle.
    $$
    The goal of this lemma is to show that $D_\Qoppa(\Lcal)$, i.e. the thick subcategory of the kernel spanned by $D_\Qoppa(z^{\gamma_S}_x)$, is a Lagrangian subcategory, making the kernel into a metabolic category. We will show that $D_\Qoppa(\Lcal)$ satisfies the hypothesis of Lemma \ref{EasierCriterionMetabolic}. \\
    
    Let us first show that $\Qoppa$ vanishes on $D_\Qoppa(\Lcal)$. Indeed $\Qoppa(D_\Qoppa(z^{\gamma_S}_x))\simeq 0$ by a direct computation, and the obstruction to $\Qoppa$ being exact on $D_\Qoppa(\Lcal)$ is precisely given by
    $$
        B_{\Qoppa}(D_\Qoppa(z^{\gamma_S}_x), D_\Qoppa(z^{\gamma_T}_y))\simeq\map(D_\Qoppa(z^{\gamma_S}_x), z^{\gamma_T}_y).
    $$
    If $S=T$, then $B_{\Qoppa}(D_\Qoppa(z^{\gamma_S}_x), D_\Qoppa(z^{\gamma_S}_y))$ is a retract of $\Qoppa(D_\Qoppa(z^{\gamma_S}_{x\oplus y}))\simeq 0$ hence vanishes itself.
    
    Using that $B_\Qoppa$ is symmetric, we are reduced to the case $T\not\subseteq S$. Let us recall that $z^{\gamma_S}_x$ is supported on those $U\subset S\cup\{n\}$ where $n\in U$ and $D_\Qoppa(z^{\gamma_T}_y)$ is supported on the edge $T\to T\cup\{n\}$ in $(\CubePunct{n})^\op$. In particular, the support of $z^{\gamma_S}_x$ and $D_\Qoppa(z^{\gamma_T}_y)$ are disjoint and furthermore we see that 
    $$
        \Map(D_\Qoppa(z^{\gamma_T}_y), z^{\gamma_S}_x)\simeq 0
    $$
    since, when computing maps in this direction, all the naturality squares have at least one of the source or target being $0$ so that they must commute for trivial reasons. In particular, we deduce:
    $$
        B_{\Qoppa}(D_\Qoppa(z^{\gamma_S}_x), D_\Qoppa(z^{\gamma_T}_y))\simeq 0.
    $$
    Since $B_\Qoppa$ is exact in both variables, we have $B_\Qoppa(D_\Qoppa(\Lcal), D_\Qoppa(\Lcal))\simeq 0$. Hence $\Qoppa$ is exact on $D_\Qoppa(\Lcal)$ and, since it vanishes on its generators, it vanishes on all of $D_\Qoppa(\Lcal)$. \\

    We now show that $\Lcal$ and $D_\Qoppa(\Lcal)$ jointly generate the kernel. Given $\gamma:i\to j$ and $\delta:j\to k$ arrows in $\CubePunct{n}$, recall that there is an exact sequence
    $$
        \begin{tikzcd}
            z^\gamma_x\arrow[r] & z^{\delta\circ\gamma}_x\arrow[r] & z^\delta_{F(\gamma)x}
        \end{tikzcd}.
    $$
    In particular, the kernel is generated by $z^\gamma_x$ for $\gamma:S\to T$ such that the cardinal $|S-T|=1$. This instance of a two-out-of-three principle means it suffices to show that the thick subcategory $\langle\Lcal,D_\Qoppa(\Lcal)\rangle$ generated by $\Lcal$ and $D_\Qoppa(\Lcal)$ contains the $z^\delta_x$ for $\delta$ of this form. In fact, we can further reduce to only edges in the source cube or between the source and the target cube; indeed, if $\eta$ is such an edge in the target cube (i.e. its source and target do not contain $n$), we can find a diagram
    $$
        
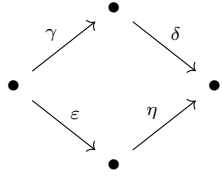
\begin{tikzcd}
            & \bullet\arrow[rd, "\delta"] \\
            \bullet\arrow[ru, "\gamma"]\arrow[rd, "\varepsilon"] & & \bullet \\
            & \bullet\arrow[ru, "\eta"]
        \end{tikzcd}
    $$
    where $\gamma$ is in the source cube and $\delta, \varepsilon$ are edges between the source and the target cube.
    
    By construction, $\Lcal$ is generated by the $z^{\gamma_S}_x$ for $\gamma_S$ a map between the source and the target cube. Let $\delta_i:\{i,n\}\to\{n\}$ be one of the map to the last vertex of the source cube. Then, the computation shows that there is an equivalence
    $$
        z^{\delta_i}_x\simeq \Sigma D_\Qoppa(z^{\gamma_{\{i\}}}_x)
    $$ 
    where $\gamma_{\{i\}}:\{i,n\}\to\{i\}$ (this is not a new notation). In particular, $z^{\delta_i}_x\in D_\Qoppa(\Lcal)$. We will now work inductively and show something a priori stronger: if $S$ is such that $n\in S$, and $T$ is such that $n\not\in T$ and $S\leq T$, we will write $\zeta^{S, T}_x$ for the section with value $F(S \leq U)(x)$ if $S\leq U\leq T$ and zero otherwise, with the canonical transition maps. In general, $\zeta^{S, T}_x$ is supported upon a cube of dimension $|S|-|T|$, and belongs to the kernel as the colimit of the induced diagram is easily checked to vanish. 
    
    We remark that all the $z^\delta_x$ for $\delta$ an edge of the source cube are of the form $\zeta^{S, T}_x$, but there are more: if $T=S-\{n\}$ we recover precisely $D_\Qoppa(z^{\gamma_S}_x)$ and, if $|S|\geq 2$, this is not of the form $z^{\delta}_x$. We now show inductively on $k=|S|-|T|$ that $\zeta^{S, T}_x\in D_\Qoppa(\Lcal)$; for $k=1$, this follows from the previous remark which identifies them with $D_\Qoppa(z^{\gamma_S}_x)$. But now generally, if $j\in S-\{n\}$, there is an exact sequence
    $$
        \begin{tikzcd}
            \zeta^{S, T\cup\{j\}}_x\arrow[r] & \zeta^{S, T}_x\arrow[r] & \zeta^{S-\{j\}, T}_x
        \end{tikzcd}
    $$
    realizing the cube spanned by $\zeta^{S, T}_x$ as the extension of two of its opposite faces. In particular, the extremal terms of the sequence are inductively known to be in $D_\Qoppa(\Lcal)$, which implies that this also holds for $\zeta^{S, T}_x$. This gives the wanted generation statement, which concludes.
\end{proof}

\begin{rmq}
    The construction of the previous Lagrangian worked regardless of the chosen direction; in particular, we have built $(n+1)$ Lagrangians in the cube.

    Note however that, if $K$ is a finite simplicial complex and $S$ is a collection of edges of the poset of faces of $K$ which are 2-by-2 disjoint and maximal in size, it need not be that the thick subcategory $\langle z^{\gamma}_x\mid\gamma\in S\rangle$ spans a Lagrangian in $\ker(\Phi^{K}_F)$. Already, in the case of $K=\CubePunct{2}$, the collection of three vertices, labeled ``$\in S$'' and colored in red, as follows:
    $$
    
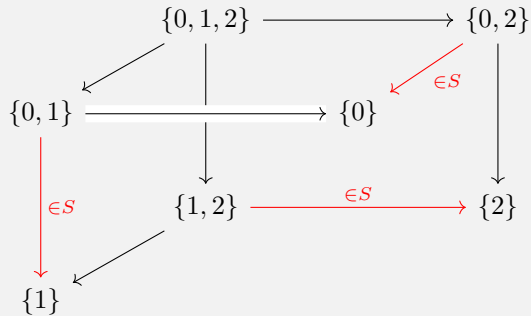
\begin{tikzcd}
    	&\{0, 1, 2\}\arrow[rr]\arrow[dl]\arrow[dd] && \{0, 2\}\arrow[dl, red, "\in S"]\arrow[dd] \\
        \{0, 1\}\arrow[dd, red, "\in S"]\arrow[rr, crossing over] && \{0\} \\
    	&\{1, 2\}\arrow[rr, red, "\in S"]\arrow[dl] && \{2\} \\
        \{1\}
    \end{tikzcd}
    $$
    does not span a Lagrangian in $\ker(\Phi^{\Delta^2}_F)$. Indeed, since all of the edges are contained in the subposet of faces of $\partial\Delta^2\subset\Delta^2$, our previous arguments imply that the thick subcategory $\langle\Lcal, D_\Qoppa(\Lcal)\rangle$ can never contain any $z^\gamma_x$ for say $\gamma:\{0, 1, 2\}\to\{0, 1\}$.
\end{rmq}

\begin{prop} \label{OplaxColimIsSentToColim}
    Let $K$ be a finite simplicial complex. For any functor $F:|K|\to\CatP$ with values in idempotent-complete categories and any $E:\CatP\to\Ecal$ Poincaré-Karoubi localizing, bordism-invariant functor, there is an equivalence
    $$
        \begin{tikzcd}
            \colim\limits_{|K|} E(F)\arrow[r, "{\simeq}"] & E(\UnP(F|_{\Face(K)^\op}))
        \end{tikzcd}
    $$
    induced by $E(F(k))\simeq E(\UnP(F|_{\Face(\{k\})^\op}))$ for $k\in|K|$.
\end{prop}
\begin{proof}
    We will work by induction on $K$ in the following way. If $K$ is a set of points, the conclusion is immediate as all such $E$ preserve direct sums and we note that there is no ambiguity to build the map. More generally, this strategy allows to reduce to the case where $K$ is connected. 
    
    There are maps $E(F(k))\simeq E(\UnP(F|_{\Face(\{k\})^\op}))\to E(\UnP(F|_{\Face(K)^\op}))$; we will show by induction that they form a cocone for $E\circ F$ and that the induced functor
    \begin{equation} \label{InducedFunctor} \tag{$\star$}
        \begin{tikzcd}
            \colim\limits_{|K|} E(F)\arrow[r] & E(\UnP(F|_{\Face(K)^\op}))
        \end{tikzcd}
    \end{equation}
    is an equivalence. \\
    
    Any finite connected simplicial complex $L$ fits into a pushout square
    $$
        \begin{tikzcd}
            K'\arrow[r]\arrow[d] & \Delta^n\arrow[d] \\
            K\arrow[r] & L
        \end{tikzcd}
    $$
    where $K'$ is a subcomplex of both $K$ and $\Delta^n$. Note that the above pushout also holds at the level of opposites of poset of faces (i.e. as categories) because simplicial complexes are glued along their faces.
    
    If $K\to L$ is an inclusion of simplicial complexes, then the induced map of opposite of the poset of (non-empty) faces of the top arrow is upwards-closed, hence by Proposition \ref{FibreCofinalInducesCartesianLift} the induced map
    $$
        \UnP(F\mid_{K})\longrightarrow\UnP(F\mid_{L})
    $$
    is a duality-preserving cartesian lift of its underlying functor of stable categories; but by Lemma \ref{UpwardsClosedInclusionAdjoint} this underlying functor is a split-Verdier inclusion given by relative left Kan extension along the inclusion of posets $K\subset L$, hence the Poincaré functor is a split Poincaré--Verdier inclusion. In particular, Lemma \ref{OplaxColimitsPreserveColimits} shows that the following square
    $$
        \begin{tikzcd}
            \UnP(F|_{K'})\arrow[r]\arrow[d] & \UnP(F\mid_{\CubePunct{n}})\arrow[d] \\
            \UnP(F\mid_{K})\arrow[r] & \UnP(F\mid_{L})
        \end{tikzcd}
    $$
    is a split Poincaré--Verdier square for $F:|L|\to\CatP$ where $\CubePunct{n}$ denotes as before the poset of faces of $\Delta^n$. Since $E$ sends such squares to pushout squares, the existence of a cocone and the fact that the induced functor \eqref{InducedFunctor} is an equivalence for $L$ now follows from the result for $\CubePunct{n}$ as well as for $K, K'$ which are simplicial complexes with strictly less faces than $L$.
    
    Hence, a straightforward induction reduces the claim to the case of $\CubePunct{n}$, $n\geq 1$. In this case, the right hand side is simply $E(F(*))$ for $F:|\CubePunct{n}|\simeq *\to\CatP$. Lemma \ref{KernelMetabolicCubes} implies that the kernel
    $$
        \ker\left(\begin{tikzcd}E\left(\UnP(F\mid_{\CubePunct{n}})\right)\arrow[r] & E\left(\colim\limits_{|\CubePunct{n}|}F\right)\end{tikzcd}\right)
    $$
    is a metabolic category. In particular, it vanishes under $E$ and we have a diagram of equivalence as follows:
    $$
        \begin{tikzcd}
            && E\left(\UnP(F\mid_{\CubePunct{n}})\right)\arrow[d, "{\simeq}"] \\
            E(F(*))\arrow[r, "{\simeq}"] & \colim\limits_{|\CubePunct{n}|}E\circ F\arrow[r, "{\simeq}"] & E\left(\colim\limits_{|\CubePunct{n}|}F\right)
        \end{tikzcd}.
    $$
    This shows both that the wanted map forms a cocone and that the induced functor is an equvialence, which concludes.
\end{proof}

For invariants that further preserve filtered colimits, Lemma \ref{GeneratedUnderFilteredColimitsByFiniteCplx} has the following immediate consequence:

\begin{cor} \label{OplaxColimIsSentToColimFinitary}
    Let $K$ be a simplicial complex. For any functor $F:|K|\to\CatP$ with values in idempotent-complete categories and any $E:\CatP\to\Ecal$ finitary Poincaré-Karoubi localizing, bordism-invariant functor, there is an equivalence
    $$
        \begin{tikzcd}
            \colim\limits_{|K|} E(F)\arrow[r, "{\simeq}"] & E(\UnP(F|_{\Face(K)^\op}))
        \end{tikzcd}
    $$
    induced by $E(F(k))\simeq E(\UnP(F|_{\Face(\{k\})^\op}))$ for $k\in|K|$.
\end{cor}

Combining this identification with the Poincaré-Verdier sequence we exhibitted earlier, we get:

\begin{thm} \label{AssemblyMapOfBIPV}
    Let $K$ be a finite simplicial complex and $F:|K|\to\CatP$ with values in idempotent-complete categories. For every Poincaré--Verdier localizing, bordism-invariant functor $E:\CatP\to\Ecal$, there is a fiber sequence
    $$
        \begin{tikzcd}
            E(\ker(\Phi^{\mathrm{P}}_F))\arrow[r] & \colim\limits_{|K|}E\circ F\arrow[r] & E(\colim\limits_{|K|}F)
        \end{tikzcd}
    $$
    where the right hand side map is the $|K|$-assembly map.
\end{thm}
\begin{proof}
    This follows from assembling the previous results: Proposition \ref{KernelStableUnderDualityFiniteCpx} shows that there is a Poincaré--Verdier sequence
    $$
        \begin{tikzcd}
            \ker(\Phi^{\mathrm{P}}_{F})\arrow[r] & \UnP(F|_{\Face(K)^\op})\arrow[r] & \UnP(F)\simeq\colim\limits_{|\Face(K)^\op|}F
        \end{tikzcd}
    $$
    whereas Proposition \ref{OplaxColimIsSentToColim} guarantees that the middle term models the left hand side of the assembly map under $E$ in such a way that the maps $F(k)\to\UnP(F|_{\Face(K)})$ induce the maps $E(F(k))\to\colim_{|K|}E\circ F$. The result follows.
\end{proof}

\begin{cor} \label{AssemblyMapOfBIPVFinitary}
    Let $K$ be a simplicial complex and $F:|K|\to\CatP$ with values in idempotent-complete categories. For every finitary Poincaré--Verdier localizing, bordism-invariant functor $E:\CatP\to\Ecal$, there is a fiber sequence
    $$
        \begin{tikzcd}
            E(\ker(\Phi^{\mathrm{P}}_F))\arrow[r] & \colim\limits_{|K|}E\circ F\arrow[r] & E(\colim\limits_{|K|}F)
        \end{tikzcd}
    $$
    where the right hand side map is the $|K|$-assembly map.
\end{cor}

There are a number of open conjectures about the assembly map for: \begin{itemize}\item quadratic L-theory, which is a Verdier-localizing bordism-invariant functor hence fits in the above picture
\item for Karoubi L-theory (which coincides with $L^{\langle -\infty\rangle}$) which is Karoubi-localizing and bordism-invariant.\end{itemize} We refer to \cite{LuckIsoConjBook} for precise statements.

For instance, the Farrell-Jones conjecture is a conjecture about the assembly map for anima of the form $BG$ for a group $G$. For quadratic L-theory and $G$ torsionfree, the Farrell-Jones conjecture ultimately reduces to the statement that
$$
    L(\ker(\Phi^{\mathrm{P}}_F))\simeq 0
$$
where $K$ is a simplicial complex whose geometric realization is $BG$. For the L-theoretic Novikov conjecture, it would suffice to prove that the map
$$
    \begin{tikzcd}
        E(\ker(\Phi^{\mathrm{P}}_F))\arrow[r] & E(\UnP(F|_K))
    \end{tikzcd}
$$
is rationally nullhomotopic. We warn the reader that we do not expect these kernels to be metabolic in general. Already for a simplicial complex modelling $S^1\vee S^1$, whose opposite of poset of faces is pictured below
$$
    
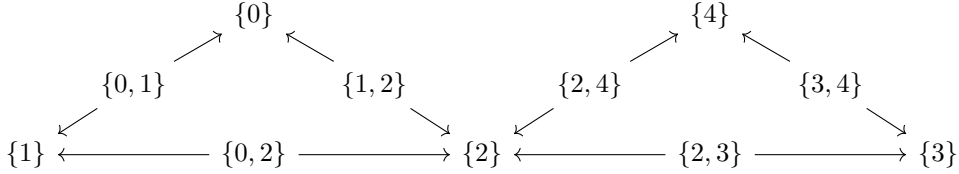
\begin{tikzcd}[row sep=small, column sep=small]
        &&\{0\} &&&& \{4\} \\
        &\{0,1\}\arrow[ru]\arrow[ld] && \{1,2\}\arrow[lu]\arrow[rd] && \{2,4\}\arrow[ru]\arrow[ld] && \{3,4\}\arrow[lu]\arrow[rd] \\
        \{1\}&& \{0,2\}\arrow[ll]\arrow[rr] && \{2\} && \{2,3\}\arrow[ll]\arrow[rr]  && \{3\} 
    \end{tikzcd},
$$
it is not possible to find a collection of edges which 2-by-2 do not have a common vertex and contain every maximal-dimensional vertex. In particular, our strategy to produce a Lagrangian subcategory fails.

However, not all hope is lost: note that metabolic categories are not closed under extensions, that is, if
$$
    \begin{tikzcd}
        (\Ccal, \Qoppa)\arrow[r] & (\Dcal, \Psi)\arrow[r] & (\Ecal, \Phi)
    \end{tikzcd}
$$
is a Verdier sequence such that the extremal terms are metabolic, then it need not be that the middle term is, precisely because being metabolic is, being in the image of, say, the right adjoint of $\fgt:\CatP\to\CatH$ and this functor is not full. Yet, $E(\Dcal, \Psi)$ will still vanish under any bordism-invariant Verdier-localizing functor $E$.

We also note that this vanishing should rely on properties of $K$, as it is not true that the assembly map in L-theory is an equivalence for every finite space, and it is only expected to hold for 1-truncated space in general.

\begin{rmq}\label{quadratic-L-theory}
    Theorem \ref{AssemblyMapOfBIPV} offers a model of the assembly map of invariants of Poincaré categories. In particular, since categories of the form $(\Ccal, \Qoppa^q_D)$ are closed under colimits in $\CatP$ and the $\Qoppa$ of $\UnH(F)$ is given by a colimit, it also provides a model for the assembly map of quadratic L-theory.

    However, the similar statement with symmetric L-theory need not hold a priori, because $(\Ccal, \Qoppa^s_D)$ are stable under limits in general, but not colimits. In specific cases however, this can hold: Lemma 4.3.8 \cite{HermKI} gives conditions for the visible symmetric structure to coincide with the symmetric one on the colimit of a sufficiently nice functor (induced by a character). Then, Proposition 6.1.9 of \loccit implies that the visible symmetric structure is indeed given by the colimit of the visible symmetric structures, at least in the case of spherical fibrations.
\end{rmq}

\section{The twisted Shaneson splitting}
\subsection{The twisted Shaneson splitting for bordism-invariant functors}

\hspace{1.2em} In this section, we deduce the twisted Shaneson splitting for bordism-invariant functors by leveraging the results of the previous section and applying them to a specific simplicial complex. Namely, let $K=\partial \Delta^2$ be the simplicial complex
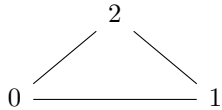
\[
\begin{tikzcd}
	& 2\arrow[rd, dash] \\
	0\arrow[ru, dash]\arrow[rr, dash] & {} & 1
\end{tikzcd}
\]
with no nondegenerate $n$-simplex for $n\geq 2$ and denote $\lambda: \catop{\Face(K)}\to |\catop{\Face(K)}|\simeq |K|\simeq S^1$ the localization at all arrows. We now show:

\begin{thm}\label{Kernel-S1-is-metabolic}
    Let $F:|K|\to \CatP$ be a functor. Then the kernel of the Poincaré functor 
    $$
        \begin{tikzcd}
            \UnH(F\circ \lambda)\arrow[r, "{\Phi^{\Ex}_{F,\lambda}}"] & \UnH(F)
        \end{tikzcd}
    $$ 
    admits a dense Poincaré subcategory which is metabolic.
\end{thm}
\begin{proof}
    Let us introduce the following notation for the poset $\catop{\Face(K)}$:
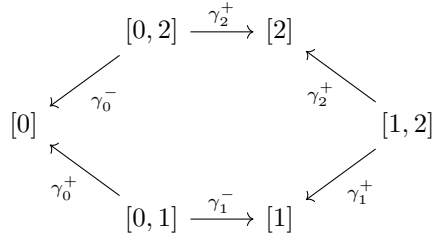
    \[\begin{tikzcd}
    	& {[0,2]} & {[2]} \\
    	{[0]} &&& {[1,2]} \\
    	& {[0,1]} & {[1]}
    	\arrow[from=1-2, to=1-3, "\gamma_2^+"]
    	\arrow[from=2-4, to=1-3, "\gamma_2^+"]
    	\arrow[from=1-2, to=2-1, "\gamma_0^-"]
    	\arrow[from=3-2, to=2-1, "\gamma_0^+"]
    	\arrow[from=2-4, to=3-3, "\gamma_1^+"]
    	\arrow[from=3-2, to=3-3, "\gamma_1^-"]
    \end{tikzcd}.\]

    Given $i=0,1,2$, we will write $i+1$ for its value mod $3$, and the same for $i-1$. This means that arrows are of the form $\gamma_i^\pm:[i, i\pm 1]\to [i]$ (up to swapping $i,i\pm 1$ inside the square brackets to match the notation).

    We know by Lemma \ref{GeneratorsKernelOfVerdierProj} that the kernel is generated as a stable category closed under retracts by objects of the form 
    $$
        z_x^{\gamma^\pm_i}:=\fib(\overline{\gamma^\pm_i})
    $$ 
    where $x\in F([i, i+1])$ and $\overline{\gamma^\pm_i}$ is the image under $\alpha$ of the cocartesian lift of $\gamma$ to $\Un(F\circ\lambda)$ with source $x$. Let us now define 
    $$
        \Lcal^\pm=\langle z_{x}^{\gamma_i^\pm}\mid i=0,1,2, x\in F([i,i+1])\rangle.
    $$ 
    It will suffice to check the conditions of Lemma \ref{EasierCriterionMetabolic} for the subcategory $\mathcal L^+\subset \ker(\Phi_{F,\lambda}^{\Ex})$ up to retracts since we only want to prove that there exists a dense subcategory which has a Lagrangian. \\

    From Remark \ref{FormulaForZGammaX} we have:
    $$
        z_x^{\gamma^+_i}(a)=
        \begin{cases*}
            0 & if $a=[i]$ \\
            x & if $a=[i,i+1]$ \\
            F(\lambda(\gamma_{i+1}^-))x & if $a=[i+1]$ \\
            0 & otherwise
        \end{cases*}
    $$
    and 
    $$
        z_x^{\gamma^-_i}(a)=
        \begin{cases*}
            F(\lambda(\gamma_{i-1}^+))x & if $a=[i-1]$\\
            x & if $a=[i-1,i]$\\
            0 & if $a=[i]$\\
            0 & otherwise.
        \end{cases*}
    $$

    %
    %
    %
    By a direct computation, we find
    $$
        D_{\Qoppa}(z^{\gamma^+_i}_x)(a)=
        \begin{cases*}
            \Sigma F(\lambda(\gamma^+_i))D_{\Qoppa_{F([i, i+1])}}x & if $a=[i]$\\
            \Sigma D_{\Qoppa_{F([i,i+1])}}x & if $a=[i,i+1]$\\
            0 & if $a=[i+1]$\\
            0 & otherwise.
        \end{cases*}
    $$
    Symmetrically, 
    $$ 
        D_{\Qoppa}(z^{\gamma^-_i}_x)(a)=
        \begin{cases*}
            0& if $a=[i-1]$\\
            \Sigma D_{\Qoppa_{F([i-1,i])}}x & if $a=[i-1,i]$\\
            \Sigma F(\lambda(\gamma_i^-))D_{\Qoppa_{F([i-1,i])}}x & if $a=[i]$\\
            0 & otherwise.
        \end{cases*}
    $$
    Hence, 
    \begin{equation}\label{duality-of-z+}
        D_\Qoppa(z^{\gamma_i^+}_x)\simeq \Sigma z^{\gamma_{i+1}^-}_{D_{\Qoppa_{F([i,i+1])}}x}
    \end{equation} 
    and symmetrically 
    $$
        D_\Qoppa(z^{\gamma_i^-}_x)\simeq \Sigma z^{\gamma^+_{i-1}}_{D_{\Qoppa_{F([i-1,i])}}x}.
    $$
    This shows that the thick subcategory $\langle \mathcal L^+, D_{\Qoppa}\mathcal L^+\rangle=\ker(\Phi^{\Ex}_{F, \lambda})$ so that we have the second part of the hypotheses of Lemma \ref{EasierCriterionMetabolic}.
    
    We now show that $\Qoppa(z_x^{\gamma_i^+})$ vanishes for $i=0$, the others case are totally symmetric; the formula is as follows:
    $$
    \Qoppa(z_x^{\gamma_0^+})\simeq\colim\left(
    \begin{tikzcd}[column sep="0.5cm"]
    	& 0 & 0 \\
    	0 &&& 0 \\
    	& \Qoppa_{F([0,1])}(x) & \Qoppa_{F([1])}(F(\lambda(\gamma_1^-))(x))
    	\arrow[from=1-2, to=1-3]
    	\arrow[from=2-4, to=1-3]
    	\arrow[from=1-2, to=2-1]
    	\arrow[from=3-2, to=2-1]
    	\arrow[from=2-4, to=3-3]
    	\arrow[from=3-2, to=3-3]
    \end{tikzcd}\right).
    $$
    Since $F$ factors through the localization at all arrows, the associated natural transformation $\Qoppa_{F(k)}\to f^*\Qoppa_{F(k')}$ is an equivalence. Therefore, the above colimit is $0$.

    Let us now show that $B_\Qoppa(z^{\gamma_i^+}_x, z^{\gamma_j^+}_y)\simeq 0$ for every $i,j=0,1,2, x\in F([i,i+1]), y\in F([j,j+1])$. If $i=j$, then $z^{\gamma_i^+}_x\oplus z^{\gamma_i^+}_y\simeq  z^{\gamma_i^+}_{x\oplus y}$ therefore $B_\Qoppa(z^{\gamma_i^+}_x, z^{\gamma_i^+}_y)$ is a retract of $\Qoppa(z^{\gamma_i^+}_{x\oplus y})\simeq 0$, hence zero itself.

    Let us now deal with the case $i\neq j$, say $i=0,j=1$; again, any other case is totally symmetric. Note that $$B_\Qoppa(z^{\gamma_0^+}_x, z^{\gamma_1^+}_y)\simeq \map(z_x^{\gamma_0^+}, D_\Qoppa(z^{\gamma_{i}^+}_y)).$$ An object of this mapping spectra is a map of diagrams in $\UnCart(F\circ \lambda)$ of the following form:

    \[\begin{tikzcd}
    	& \ 0 & \ 0 \\
    	0 &&& 0 \\
    	& x \ & {F(\lambda(\gamma_{1}^+))x} \\
    	& \ 0 & \ 0 \\
    	0 &&& {\Sigma D_{\Qoppa_{F([1,2])}}y} \\
    	& 0 \ & {\Sigma F(\lambda(\gamma_1^+))D_{\Qoppa_{F([1,2])}}y}
    	\arrow[from=1-2, to=1-3]
    	\arrow[from=1-2, to=2-1]
    	\arrow[shift left=2, from=1-2, to=3-2, dash]
    	\arrow[shift left=2, from=1-3, to=3-3, dash]
    	\arrow[from=2-1, to=5-1]
    	\arrow[from=2-4, to=1-3]
    	\arrow[from=2-4, to=3-3]
    	\arrow[from=2-4, to=5-4]
    	\arrow[from=3-2, to=2-1]
    	\arrow[from=3-2, to=3-3]
    	\arrow[shift right=2, from=3-2, to=6-2, crossing over]
    	\arrow[shift left=2, from=3-2, to=4-2]
    	\arrow[shift right=2, from=3-3, to=6-3, crossing over]
    	\arrow[shift left=2, from=3-3, to=4-3]
    	\arrow[from=4-2, to=4-3]
    	\arrow[from=4-2, to=5-1]
    	\arrow[from=5-4, to=4-3]
    	\arrow[from=5-4, to=6-3]
    	\arrow[from=6-2, to=5-1]
    	\arrow[from=6-2, to=6-3]
    \end{tikzcd}.\]
    The commutativity of the front square induces a factorization of the map 
    $$
        \begin{tikzcd}
            F(\lambda(\gamma_{1}^+))x\arrow[r] & \Sigma F(\lambda(\gamma_1^+))D_{\Qoppa_{F([1,2])}}y
        \end{tikzcd}
    $$
    through the cokernel of the top map $x\to F(\lambda(\gamma_{1}^+))x$, which is an equivalence hence has vanishing cokernel. Therefore the wanted mapping spectrum is null.

    This shows that $\Qoppa$ is exact on $\mathcal L^+$ and, since it vanishes on the generators $z_x^{\gamma_i^+}$ it is zero on $\mathcal L$. In particular, we have verified the hypotheses of Lemma \ref{EasierCriterionMetabolic} and hence the proof of Theorem \ref{Kernel-S1-is-metabolic} is complete.
\end{proof}

\begin{cor} \label{AssemblyBZofBIPVisEquivalence}
    Let $F:\tB\Z\to\CatP$ be a functor such that $F(*)$ is idempotent-complete, and $E:\CatP\to\Ecal$ be a bordism-invariant Poincaré--Verdier localizing functor. Then, the assembly map
    $$
        \begin{tikzcd}
            \colim\limits_{\tB\Z}E\circ F\arrow[r] &  E(\colim\limits_{\tB\Z}F)
        \end{tikzcd}
    $$ 
    is an equivalence in $\Ecal$.
\end{cor}

\begin{proof}
    This follows from combining Theorem \ref{AssemblyMapOfBIPV} and Theorem \ref{Kernel-S1-is-metabolic}.
\end{proof}

In \cite[Theorem 5.1]{Ranicki-ALTIII}, the author considers a ring with involution $(R, \tau)$ equipped with an automorphism $\alpha$ as well as $T\subset \widetilde{\textup{K}}_0(R)$ a $\tau$-invariant subgroup. He proves that there exists a long exact sequence of abelian groups
$$
    \begin{tikzcd}[column sep=0.4cm]
        \dots\arrow[r] & \tL^{\tq,(1-\alpha)^{-1}T}_n(R, \tau)\arrow[r, "{1-\alpha}"] & \tL^{\tq,T}_n(R, \tau)\arrow[r, "\bar\varepsilon"] & \tL^{\tq,\bar\varepsilon T}_n(R_\alpha[t^{\pm}], \tau^-)\arrow[r] & \tL^{\tq,(1-\alpha)^{-1}T}_{n-1}(R, \tau)\arrow[r]&\dots .
    \end{tikzcd}
$$
where $\tau^-$ denotes the involution $\tau$ extended to Laurent polynomials via $t\mapsto t^{-1}$. 

We prove a version of this result for Karoubi L-theory $\bL:\CatP\to\Sp$, the initial Karoubi-localizing invariant under L-theory, as described in \cite{HermKIV}. In particular, this yields a twisted version of \cite[Example 2.1.14]{HermKIV}.

\begin{cor}\label{Shaneson-Karoubi}
    Let $(R, \tau)$ be a ring with involution, equipped with an automorphism $\alpha$. Let $(R_\alpha[t^\pm], \tau^-)$ be the ring of $\alpha$-twisted Laurent polynomials equipped with the involution sending $t\mapsto t^{-1}$. Then for every $m\in [-\infty, \infty]$ there exists a long exact sequence of abelian groups
    $$
        \begin{tikzcd}[column sep=0.6cm]
            \dots\arrow[r] & \bL^{\geq m}_n(R)\arrow[r, "{1-\alpha}"] & \bL^{\geq m}_n(R, \tau)\arrow[r] & \bL^{\geq m}_n(R_\alpha[t^{\pm}], \tau^-)\arrow[r] & \bL^{\geq m}_{n-1}(R, \tau)\arrow[r]&\dots
        \end{tikzcd}
    $$
    where $\bL$ denotes Karoubi L-theory.
\end{cor}
\begin{proof}
    Let $F:\tB\Z\to \CatP$ be the functor with value $(\mathcal D^f(R), \Qoppa_\tau^{\geq m})$ at the point and sending the generator of $\tB\Z$ to the automorphism induced by $\alpha$, which we still denote by $\alpha$. By \cite[§4.3]{HermKI} (see also \cite[Remark 2.2.1]{HermKIV}), we can identify up to idempotent-completion the colimit of $F$ as follows:
    \begin{equation}
        \Idem(\colim_{\tB\Z} F)\simeq (\Perf(R_\alpha[t^{\pm}]), \Qoppa_{\tau^-}^{\geq m}).
    \end{equation}
    Hence by invariance under taking idempotent completions of Karoubi L-theory, we have 
    $$
        \bL(\colim_{\tB\Z} F)\simeq \bL^{\geq m}_n(R_\alpha[t^{\pm}], \tau^-)
    $$
    where the right hand side is simply our notation for $\bL(\Perf(R_\alpha[t^{\pm}]), \Qoppa_{\tau^-}^{\geq m})$. By Corollary \ref{AssemblyBZofBIPVisEquivalence} applied to $\Idem(F), E=\bL$, we thus obtain 
    $$
        \bL(\Perf(R_{\alpha}[t^\pm], \Qoppa^{\geq m}_{\tau^-}))\simeq \colim_{\tB\Z}\bL\circ F.
    $$ 
    Finally, we can compute this colimit as
    \begin{gather*}
            \colim_{\tB\Z}\bL\circ F\simeq\colim \left(\begin{tikzcd}[ampersand replacement=\&]\bL^{\geq m}(R, \tau) \& \arrow[l, "{\textup{pr}_1}"']\bL^{\geq m}(R, \tau)\oplus\bL^{\geq m}(R, \tau)\arrow[r, "{\alpha\circ \textup{pr}_2}"] \& \bL^{\geq m}(R, \tau)\end{tikzcd}\right) \\
            \simeq\cofib \left(\begin{tikzcd}[ampersand replacement=\&]\bL^{\geq m}(R,\tau)\arrow[r, "{\id-\alpha}"] \& \bL^{\geq m}(R,\tau)\end{tikzcd}\right).
    \end{gather*}
    Taking homotopy groups leads to the claimed long exact sequence.
\end{proof}

\subsection{Towards the general twisted Shaneson splitting}

\hspace{1.2em} Let $\Ical$ be the following category
\begin{equation}\label{DiagramIcal}
    \begin{tikzcd}
        \bullet\arrow[loop left, "(-1)", phantom]\arrow[r, bend left=30, "\alpha_+"]\arrow[r, bend right=30, "\alpha_-"'] & \bullet \arrow[loop right, "(+1)", phantom]
    \end{tikzcd}.
\end{equation}
Recall that in our convention, $\TwAr(\Delta^1)$ is the category depicted as follows:
$$
    \begin{tikzcd}
        (0\to 0) & (0\to 1)\arrow[l]\arrow[r]& (1\to 1)
    \end{tikzcd}
$$
Hence, there is a pushout square of categories
$$
    \begin{tikzcd}
        \Delta^0\coprod\Delta^0\arrow[r]\arrow[d] & \TwAr(\Delta^1)\arrow[d] \\
        \Delta^0\arrow[r] & \Ical
    \end{tikzcd}
$$  
where the bottom map is the inclusion of $(+1)$ in the description (\ref{DiagramIcal}). We deduce the following:

\begin{prop} \label{PoincareSquareUnPIcal}
    Let $F:\Ical\to\CatP$ be a functor. There is a pushout square in $\CatP$:
    $$
        \begin{tikzcd}
            F(1)\oplus F(1)\arrow[r]\arrow[d] & \UnP(F\mid_{\TwAr(\Delta^1)})\arrow[d, "p"] \\
            F(1)\arrow[r, "i_{1}"] & \UnP(F)
        \end{tikzcd}
    $$
    whose horizontal arrows are split Poincaré--Verdier inclusions. In particular, $\UnP(F)$ is Poincaré.
\end{prop}
\begin{proof}
    The existence of the above square and the fact that it is a pushout of hermitian categories follow from Lemma \ref{OplaxColimitsPreserveColimits}.

    Since we already know the square is a pushout in $\CatH$, the closure of $\CatP$ under colimits in $\CatH$ reduces the problem to checking whether the left vertical map and the top horizontal map are Poincaré functors between Poincaré categories. The hermitian categories $F(1)$ and $F(1)\oplus F(1)$ are obviously Poincaré and Lemma \ref{UnPPoincaréCubePunctN} guarantees this also holds for $\UnP(F\mid_{\TwAr(\Delta^1)})$ since $\TwAr(\Delta^1)\simeq\CubePunct{1}$.
    \begin{itemize}
        \item[(i)] The map $F(1)\oplus F(1)\to F(1)$ is the fold map, hence Poincaré for reasons unrelated to $F$.
        \item[(ii)] The functor $\Delta^0\coprod\Delta^0\to\TwAr(\Delta^1)$ satisfies the hypotheses of Proposition \ref{FibreCofinalInducesCartesianLift}, hence $F(1)\oplus F(1)\to\UnP(F\mid_{\TwAr(\Delta^1)})$ is a duality-preserving cartesian lift of the underlying functor of exact categories. To show that it is a split-Poincaré--Verdier inclusion, it suffices to show that
        $$
            F(1)\oplus F(1)\to\StabUn(F\mid_{\TwAr(\Delta^1)})
        $$
        is fully-faithful with a right adjoint. We already have the adjunction statement by Lemma \ref{AlphaiasAdjoints} and the fully-faithfulness can be checked on the formula using $\Map_{\TwAr(\Delta^1)}(i, i)\simeq *$ for every $i\in\TwAr(\Delta^1)$.
    \end{itemize}
    Since $\CatP\to\CatH$ reflects colimits, the square is a pushout square in $\CatP$. This implies that the bottom horizontal functor is a Poincaré--Verdier inclusion, which concludes.
\end{proof}

\begin{prop} \label{CommonValueKernelUnPIcal}
    Let $F:\Ical\to\CatP$ be a functor. The common value of the cokernel of the Poincaré--Verdier inclusions of Proposition \ref{PoincareSquareUnPIcal} is $(F(-1), \Sigma\Qoppa_{F(-1)})$. Differently stated, the Poincaré--Verdier inclusion $i_1:F(1)\to\UnP(F)$ fits into a split Poincaré--Verdier sequence
    $$
        \begin{tikzcd}
            F(1)\arrow[r, "i_{1}"] & \UnP(F)\arrow[r, "\ev_{-1}"] & (F(-1), \Sigma\Qoppa_{F(-1)})
        \end{tikzcd}
    $$
\end{prop}
\begin{proof}
    By Proposition \ref{PoincareSquareUnPIcal}, we know that $i_1$ is a split Poincaré--Verdier inclusion whose cokernel is the same as $F(1)\oplus F(1)\to\UnP(F\mid_{\TwAr(\Delta^1)})$. We are thus left with the task of identifying the latter. We claim that
    $$
        \ev_{-1}:\UnP(F\mid_{\TwAr(\Delta^1)})\longrightarrow(F(-1), \Sigma\Qoppa_{F(-1)})
    $$
    is the quotient functor (of the latter map, but consequently also of the former). First let us remark that $\ev_{-1}$ is indeed a Poincaré functor thanks to Remark \ref{AlphaNotPoincareButEvYes}.

    Using that $\TwAr(\Delta^1)$ is obtained by gluing two copies of $\Delta^1$ along three copies of $\Delta^0$, the explicit description of \ref{ExplicitStabUnStronglyFinite} implies that $\StabUn(F\mid_{\TwAr(\Delta^1)})$ fits in the following pullback square
    $$
        \begin{tikzcd}[column sep=huge]
            \StabUn(F\mid_{\TwAr(\Delta^1)})\arrow[d, "{(\ev_{1}, \ev_{1}, \ev_{-1})}"]\arrow[rr] && \Ar(F(1))^2\arrow[d, "{(s, t, s, t)}"] \\
            F(1)\oplus F(1)\oplus F(-1)\arrow[rr, "{(\id, \id, f(\alpha^+), f(\alpha^-))}"]& & F(1)^4
        \end{tikzcd}.
    $$
    In particular, the underlying functor of stable categories of $\ev_{-1}$ is indeed the cokernel of $F(1)\oplus F(1)\to \StabUn(F\mid_{\TwAr(\Delta^1)})$. 
    
    To conclude, we have to check that
    $$
        (\ev_{-1}^\op)_!\Qoppa_F\longrightarrow\Sigma\Qoppa_{F(-1)}
    $$
    is an equivalence. Since the formula for $(\ev_{-1}^\op)_!\Qoppa_F$ is the colimit over $\TwAr(\Delta^1)$ --- i.e. the pushout --- of the following diagram:
    $$
        \begin{tikzcd}
            & (\ev_{-1}^\op)_!(\ev_{1}^\op)^*\Qoppa_{F(1)} \\
            (\ev_{-1}^\op)_!(\ev_{-1}^\op)^*\Qoppa_{F(-1)}\arrow[ru]\arrow[rd] \\
            & (\ev_{-1}^\op)_!(\ev_{1}^\op)^*\Qoppa_{F(1)}
        \end{tikzcd}
    $$
    it suffices to check that $(\ev_{-1}^\op)_!(\ev_{-1}^\op)^*\simeq\id$ and $(\ev_{-1}^\op)_!(\ev_{1}^\op)^*\simeq 0$. The first statement follows from the fully-faithfulness of the right adjoint of $\ev_{-1}$, itself a consequence of the formula of Lemma \ref{AlphaiasAdjoints} and of the fact that every point of $\TwAr(\Delta^1)$ has no non-trivial endomorphisms. Since $(\ev_{1}^\op)^* \simeq(\alpha_1^\op)_!$, the second statement follows from checking that $\ev_{-1}\circ\;\alpha_{1}\simeq0$ which we can also see in the formula of Lemma \ref{AlphaiasAdjoints}, as $\Map_{\TwAr(\Delta^1)}(1, -1)=\emptyset$. This concludes.
\end{proof}

Remark that $\TwAr(\Delta^1)$ is contractible and the map $\TwAr(\Delta^1)\to\Ical\to \tB\Z$ therefore factorizes through $*\to B \Z$. Hence, given a functor $F:\tB\Z\to\CatP$, we can form the comparison functor between the oplax colimit and the colimit of two different functors which gives rises to:

\begin{lmm}
    Let $F:\tB\Z\to\CatP$ be a functor. There is a diagram with Poincaré--Verdier sequences as rows as follows:
    $$
        \begin{tikzcd}
            \ker(\Phi^{\mathrm{P}}_{F\mid_{\TwAr(\Delta^1)}})\arrow[r]\arrow[d] & \UnP(F\mid_{\TwAr(\Delta^1)})\arrow[r]\arrow[d] & F(*)\arrow[d] \\
            \ker(\Phi^{\mathrm{P}}_{F})\arrow[r]& \UnP(F\mid_{\Ical})\arrow[r] & \colim\limits_{\tB\Z} F
        \end{tikzcd}.
    $$
    Moreover, the first kernel $\ker(\Phi^{\mathrm{P}}_{F\mid_{\TwAr(\Delta^1)}})$ is a metabolic category.
\end{lmm}
\begin{proof}
    The functor $\TwAr(\Delta^1)\to\Delta^0\simeq *$ is the localisation at all the arrows of a poset of faces of a simplicial complex, hence the top sequence is a Poincaré--Verdier sequence by an instance of Proposition \ref{KernelStableUnderDualityFiniteCpx}. We postpone the computation that shows the stability of the kernel $\ker(\Phi^{\mathrm{P}}_{F})$ to a latter section (see Lemma \ref{KernelStableUnderDuality}), for the only reason of gathering all of the computations about it in a single place (this will imply that the bottom sequence is a Poincaré--Verdier sequence as well). 
    
    The compatibility of the sequences is ensured by the functoriality of oplax colimits, see Proposition \ref{OplaxColimitsPreserveColimits}. Finally, the metabolicity of the top kernel is a direct consequence of Lemma \ref{KernelMetabolicCubes}, since $\TwAr(\Delta^1)$ is the opposite of the poset of faces of the simplicial complex $\begin{tikzcd}\bullet_0\arrow[r, dash] &\bullet_1\end{tikzcd}$.
\end{proof}

We can combine the square of Proposition \ref{PoincareSquareUnPIcal} and the cokernel of \ref{CommonValueKernelUnPIcal} with the rectangle of the above lemma to get the following:
$$
    \begin{tikzcd}[row sep=small, column sep=small]
        0\arrow[rr]\arrow[dd]\arrow[rd] && F(*)\oplus F(*)\arrow[rr]\arrow[dd]\arrow[rd] && F(*)\oplus F(*)\arrow[dd]\arrow[rd] \\
        &\ker(\Phi^{\mathrm{P}}_{F\mid_{\TwAr(\Delta^1)}})\arrow[rr, crossing over] && \UnP(F\mid_{\TwAr(\Delta^1)})\arrow[rr, crossing over]&& F(*)\arrow[dd] \\
        0\arrow[rr]\arrow[rd] && F(*)\arrow[rr]\arrow[rd] && F(*)\arrow[rd] \\
        &\ker(\Phi^{\mathrm{P}}_{F})\arrow[rr]\arrow[from=uu, crossing over] && \UnP(F\mid_{\Ical})\arrow[rr]\arrow[from=uu, crossing over] && \colim\limits_{\tB\Z} F
    \end{tikzcd}.
$$
The rightmost face is cocartesian by abstract nonsense, and every horizontal sequence is a Poincaré--Verdier sequence. Recall from Lemma \ref{PoincareSquareUnPIcal} that the middle slanted square is a pushout in $\CatP$, whose horizontal arrows are Poincaré--Verdier inclusions.

Thanks to these observations, for any Poincaré--Verdier localizing $E$ there is a diagram of short-exact sequences
\begin{equation}\label{eq-ComparisonReducesToKernels}
    \begin{tikzcd}
        E(\ker(\Phi^{\mathrm{P}}_{F\mid_{\TwAr(\Delta^1)}}))\arrow[r]\arrow[d] & E(\UnP(F\mid_{\Ical}))\arrow[r]\arrow[d, equal]& \colim_{\tB\Z}E\circ F\arrow[d]\\
        E(\ker(\Phi^{\mathrm{P}}_{F}))\arrow[r] & E(\UnP(F\mid_{\Ical}))\arrow[r] & E(\colim_{\tB\Z}F)
    \end{tikzcd},
\end{equation}
where the rightmost map is the colimit-comparison map, that is the assembly map for $\tB\Z$. In other words, the failure for the $\tB\Z$-assembly map to be an equivalence is completely controlled by the failure of the map between the kernels to be an equivalence.

The goal of the next section is to study the Poincaré functor:
$$
    \begin{tikzcd}
        \ker(\Phi^{\mathrm{P}}_{F\mid_{\TwAr(\Delta^1)}})\arrow[r, "{\Psi}"] & \ker(\Phi^{\mathrm{P}}_{F}).
    \end{tikzcd}
$$ 
We already know that $\ker(\Phi^{\mathrm{P}}_{F\mid_{\TwAr(\Delta^1)}})$ is hyperbolic, and we will show that the whole map is in fact in the image of $\Hyp$, and that it admits a retraction which becomes compatible with the previous diagram after applying a Poincaré--Verdier localizing invariant.

\subsection{The nil-term(s) and the twisted Shaneson splitting}

\hspace{1.2em} We want to show that $\ker(\Phi^{\mathrm{P}}_{F})$ is metabolic (with the same notation as in the previous section). From Lemma \ref{UnPNonDegenerateAndFormula}, we have that for $x\in F(-1)$:
\begin{align*}
    \alpha_{-1}(x) & \simeq \begin{tikzcd}[ampersand replacement = \&]\Big[(-1,\lim_{\alpha\in \Ical_{-1//-1}}^{F(-1)}x) \& (1,\lim_{\alpha\in \Ical_{-1//1}}^{F(1)}F(\alpha)x)\arrow[l,shift right]\arrow[l,shift left] \Big]\end{tikzcd} \\
                   & \simeq \begin{tikzcd}[ampersand replacement = \&]\Big[(-1,x) \& (1,F(\alpha^+)x\oplus F(\alpha^-)x)\arrow[l,shift right]\arrow[l,shift left] \Big]\end{tikzcd}
\end{align*}
with the first map (resp. the second map) given by the first (resp. second) projection of the direct sum followed by the cartesian lift of $\alpha^+$ (resp. $\alpha^-$) at $x\in F(-1)\subset\UnCart(F)$. Here, for simplicity, we chose to write a section $\Ical^\op\to\UnCart(F)$ by the diagram it represents. Analogously, one can compute:
$$
    \alpha_1(F(\alpha^+)x)\simeq\begin{tikzcd}\Big[(-1,0) & (1,F(\alpha^+)x)\arrow[l,shift right]\arrow[l,shift left] \Big]\end{tikzcd} $$
and
$$
    \alpha_1(F(\alpha^-)x)\simeq\begin{tikzcd}\Big[(-1,0) & (1,F(\alpha^-)x)\arrow[l,shift right]\arrow[l,shift left] \Big]\end{tikzcd}
$$
where the two maps are the unique such maps over either $\alpha^+$ or $\alpha^-$. The oplax structure of $\StabUn(F)$ means there are natural maps
$$
    \begin{tikzcd}
        \alpha_{-1}(x)\arrow[r] & \alpha_1(F(\alpha^\pm)x)
    \end{tikzcd}
$$
whose fiber is given by:
$$
    z^{\pm}_x\coloneqq\fib(\begin{tikzcd}\alpha_{-1}(x)\arrow[r] &  \alpha_1(F(\alpha^\pm)x)\end{tikzcd})=\begin{tikzcd}\Big [(-1,x)& (1,F(\alpha^\mp)x)\arrow[l,shift right]\arrow[l,shift left]\Big ].\end{tikzcd}
$$
Here, in the case of $z^+_x$, the first map is the cartesian lift of $\alpha^-$ and the second is zero. In the case of $z^-_x$, the first map is zero and the second is the identity. 

\begin{lmm} \label{KernelStableUnderDuality}
    The kernel $\ker(\Phi^{\Ex}_{F})$ is generated under finite colimits and retracts by the $z^{\pm}_x$ for $x\in F(-1)$. Moreover, we have 
    $$
        D_{\Qoppa_F}(z^\pm_x)\simeq \Sigma z^{\mp}_{D_{F(-1)}x}.
    $$
    In particular, the kernel is stable under the duality of $\StabUn(F)$. 
\end{lmm}
\begin{proof}
    The generation statement follows directly from Lemma \ref{GeneratorsKernelOfVerdierProj}. Using Proposition \ref{UnPNonDegenerateAndFormula}, we can compute
    $$
        D_{\Qoppa_F}(z^\pm_x)(-1)\simeq\colim^{F(-1)}\Big[\begin{tikzcd} D_{F(-1)}x \arrow[r,shift left]\arrow[r, shift right] & 0\Big]\end{tikzcd})\simeq\Sigma D_{F(-1)}x
    $$
    where the superscript $F(-1)$ indicates we take the colimit in the stable category $F(-1)$, and 
    \begin{align*}
        D_{\Qoppa_F}(z^\pm_x)(1) & \simeq \colim^{F(1)}\Big[\begin{tikzcd}[ampersand replacement = \&] D_{F(1)}(F(\alpha^+)x\oplus F(\alpha^-)x) \arrow[r,shift left]\arrow[r, shift right] \& D_{F(1)}F(\alpha^\mp)x\Big]\end{tikzcd} \\
        &\simeq\Sigma D_{F(1)}(F(\alpha^\pm)x) \\
        &\simeq\Sigma F(\alpha^\pm)(D_{F(-1)}x)
    \end{align*}
    where the last equality follows from the fact that $F(\alpha^{\pm})$ is duality-preserving. Comparing with the formula, we get the wanted formula on the duality, and the last statement is straightforward.
\end{proof}

From the above lemma, it appears that there are two flavours of objects in the kernel mentioned in the statement, which warrants the following definition:

\begin{defi}
    Let $\Ncal^+$ denote the full stable subcategory of $\StabUn(F)$ which is closed under retracts generated by the $z^{+}_x$ for $x\in F(-1)$ (respectively $\Ncal^-$ for $z^{-}_x$ still with $x\in F(-1)$).
\end{defi}

\begin{lmm} \label{QuadraticVanishing}
    The quadratic functor on $\ker\Phi^{\mathrm{p}}_{F}$ vanishes on both $\Ncal^+$ and $\Ncal^-$. Moreover, the duality induces an equivalence $\catop{(\Ncal^+)}\simeq\Ncal^-$.
\end{lmm}
\begin{proof}
    The second half of the statement follows from Lemma \ref{KernelStableUnderDuality} since the duality is an anti-involution of the kernel. To check the first half, first recall from \ref{UnPNonDegenerateAndFormula}:
    \begin{equation}\label{coequalizer}
    \Qoppa_{F}(z^\pm_x)\simeq\colim\begin{tikzcd}\Big[\Qoppa_{F(-1)}(x)\arrow[r,shift left, "\eta_x"]\arrow[r,shift right, "0"']&\Qoppa_{F(1)}(F(\alpha^\mp)x)\Big]\end{tikzcd}
    \end{equation}
    where $\eta^\mp_x$ is the value at $x$ of the natural transformation $\eta^\mp:\Qoppa_{F(-1)}\to \Qoppa_{F(1)}\circ \catop{F(\alpha^\mp)}$ provided by Lemma \ref{DatumOfLift}. Since it factors through $\tB\Z$, $F$ sends all arrows of $\Ical$ to equivalences (this is actually the only point we use that property) hence, in particular, the two possible natural transformations $$\Qoppa_{F(-1)}\to \Qoppa_{F(1)}\circ \catop{F(\alpha^\mp)}$$ are equivalences.
    Therefore, the top arrow in the diagram \eqref{coequalizer} is an equivalence, and hence the coequalizer is zero. 
    
    We have checked that $\Qoppa$ vanishes on the generators of $\Ncal^+, \Ncal^-$ but $\Qoppa$ is only quadratic and not exact so this does not suffices. However, the only obstruction to $\Qoppa$ being exact is $B_\Qoppa$, which is linear in each variable and thus to conclude, it suffices to check that $B_\Qoppa(z^+_x, z^+_y)\simeq 0$. Now, $B_\Qoppa(z^+_x, z^+_y)$ is a retract of $\Qoppa(z^+_x\oplus z^+_y)$; but since $z^+_x\oplus z^+_y\simeq z^+_{x\oplus y}$, the previous computation shows that $\Qoppa(z^+_x\oplus z^+_y)\simeq 0$, which concludes.
\end{proof}

\begin{lmm} \label{OrthogonalComplement}
    The two subcategories $\Ncal^+$ and $\Ncal^-$ are in orthogonal complement inside $\ker\Phi^{\Ex}_{F}$.
\end{lmm}
\begin{proof}
    Let us simply write $\Qoppa$ for the quadratic functor on $\ker\Phi^{\mathrm{p}}_{F}$. Remark that
    $$
        \map(z^+_x, z^-y)\simeq\map(z^+_x, D_\Qoppa(z^+_y))\simeq B_\Qoppa(z^+_x, z^+_y).
    $$
    As previously, $B_\Qoppa(z^+_x, z^+_y)$ is a retract of $\Qoppa(z^+_{x\oplus y})$, which we have shown to vanish in the previous lemma, hence
    $$
        \map(z^+_x, z^-_y)\simeq 0.
    $$
    Since the duality is its own adjoint, the above computation also shows that
    $$
        \map(z^-_x, z^+_y)\simeq 0
    $$
    which concludes.
\end{proof}

We have now all that we need to show the following:

\begin{prop} \label{KernelHyperbolic}
    Let $F:\tB\Z\to\CatP$. There is a Poincaré--Verdier functor, whose underlying exact functor is a dense inclusion
    $$
        \begin{tikzcd}
            \Hyp(\Ncal^+)\arrow[r] & \ker(\Phi^{\mathrm{p}}_{F})
        \end{tikzcd}.
    $$
    This is an equivalence if the composite $F:\tB\Z\to\CatEx$ has values in idempotent-complete stable categories.
\end{prop}
\begin{proof}
    By Lemma \ref{OrthogonalComplement}, the functor $\Ncal^+\oplus\Ncal^-\to\ker\Phi^{\Ex}_{F, i}$ is fully-faithful. Moreover, by Lemma \ref{QuadraticVanishing}, the quadratic functor on $\ker(\Phi^{\mathrm{p}}_{F, i})$ vanishes on both summands so its restriction to the direct sum must be hyperbolic. Lemma \ref{GeneratorsKernelOfVerdierProj} guarantees that the inclusion $\Ncal^+\oplus\Ncal^-\to\ker\Phi^{\Ex}_{F, i}$ has dense image, which gives the first part. \\
    
    Suppose now that the underlying stable categories in which $F$ takes values are idempotent-complete, then so is $\StabUn(F)$ by Proposition \ref{IdempotentCompleteStabUn}; note that this applies both to $F$ and $F\circ i$ so that the kernel $\ker(\Phi^{\Ex}_{F, i})$ is also idempotent-complete, as such categories are closed under limits in $\CatEx$. Since they are closed under retracts in the kernel, $\Ncal^+$ and $\Ncal^-$ are idempotent-complete as well and thus so is their sum $\Ncal^+\oplus\Ncal^-$. Consequently, the above map is an equivalence as wanted.
\end{proof}

For the following, we will need the following folklore result, related to results of Schwede--Shipley:

\begin{lmm} \label{MapsFromStableCatGeneratedByOneElement}
    Let $\Ccal$ be a small stable category generated under finite colimits and retracts by a single object $G\in\Ccal$, and denote $B\End_\Ccal(G)$ the one-point category with $\End_\Ccal(G)$ as endomorphisms. The functor $\alpha:B\End_\Ccal(G)\to\Ccal$ induces an equivalence
    $$
        \begin{tikzcd}
            \alpha^*:\FunEx(\Ccal, \Dcal)\arrow[r, "\simeq"] & \Fun(B\End_\Ccal(G), \Dcal)
        \end{tikzcd}
    $$
    for every idempotent-complete $\Dcal$.
\end{lmm} 
\begin{proof}
    Under these hypotheses, $\Ind(\Ccal)$ is generated under colimits by $G$ and therefore
    $$
        \begin{tikzcd}
            \alpha^*:\FunL(\Ind(\Ccal), \Dcal)\arrow[r] & \Fun(B\End_\Ccal(G), \Dcal)
        \end{tikzcd}
    $$
    is conservative for any presentable stable $\Dcal$. A priori, this functor has a left adjoint given by left Kan extension followed by forcing the functor to commute with small colimits.
    
    We now show that the left Kan extension $\alpha_!$ already lands in $\FunL(\Ind(\Ccal), \Dcal)$ so that this second operation is not necessary; for this we adapt a criterion of Rezk, see \cite{MO-Rezk}, to the context of stable categories. Write $\Pcal^{st}$ for the left adjoint $\Cat\to\CatEx$ to the inclusion and let us shorten $\Ccal_0:=B\End_\Ccal(G)$. Then, there is a unique colimit-preserving $\overline{\alpha}:\Ind(\Pcal^{st}(\Ccal_0))\to\Ccal$ whose restriction to $\Ccal_0$ is $\alpha$. In particular, there is a commutative diagram
    $$
        \begin{tikzcd}
            \FunL(\Ind(\Ccal), \Dcal)\arrow[r, "{\alpha^*}"]\arrow[rd, "{\overline{\alpha}^*}"] & \Fun(\Ccal_0, \Dcal) \\
            & \FunL(\Ind(\Pcal^{st}(\Ccal_0)), \Dcal)\arrow[u, "\simeq"]
        \end{tikzcd}.
    $$
    Hence, passing to left adjoints, for any $g:B\End_\Ccal(G)\to\Dcal$, we see that $\alpha_!g$ is equivalently given as the composite $\Ind(\Pcal^{st}(g))\circ R_\alpha$ where $\Ind(\Pcal^{st}(g))$ is the unique colimit-extension of $g$ to $\Ind(\Pcal^{st}(\Ccal_0))$ and $R_\alpha$ is the following functor:
    \begin{align*}
        \Ind(\Ccal)\longrightarrow & \Ind(\Pcal^{st}(\Ccal_0))\simeq\FunEx(\Pcal^{st}(\Ccal_0)^\op, \Sp) \\
        X\longmapsto & \map_{\Ind(\Ccal)}(\alpha(-), X)
    \end{align*}
    formally induced by taking the right adjoint of the colimit-preserving extension of $\alpha$ to $\Ind(\Pcal^{st}(\Ccal_0))$. In particular, it suffices to check that $R_\alpha$ is colimit-preserving for all of the $\alpha_!g$ to be, and this follows directly from the fact that $\map_{\Ind(\Ccal)}(\alpha(C), -)$ commutes with colimits for every $C\in\Pcal^{st}(\Ccal_0)$ because they are compact in $\Ind(\Ccal)$ since our original $\alpha:B\End_\Ccal(G)\to\Ccal$ lands in $\Ccal\subset\Ind(\Ccal)^\omega$. \\

    In particular, the above showed that, for every presentable stable $\Ecal$, there is an equivalence
    $$
        \begin{tikzcd}
            \alpha^*:\FunEx(\Ccal, \Ecal)\arrow[r, "\simeq"] & \Fun(B\End_\Ccal(G), \Ecal)
        \end{tikzcd}.
    $$
    Now given a small stable $\Dcal$, instantiating the above statement with $\Ecal=\Ind(\Dcal)$, we get that
    $$
        \begin{tikzcd}
            \alpha^*:\FunEx(\Ccal, \Dcal)\arrow[r] & \Fun(B\End_\Ccal(G), \Dcal)
        \end{tikzcd}
    $$
    is fully-faithful. But note that, given $f:B\End_\Ccal(G)\to\Dcal$, for every $X\in\Ccal$ such that $X\simeq\colim_I G$ with $I$ finite, the unique extension $\overline{f}:B\End_\Ccal(G)\to\Ind(\Dcal)$ must satisfy
    $$
        \overline{f}(X)\simeq\colim_{I}f(G)
    $$
    hence $\overline{f}(X)\in\Ccal$. Since every object of $\Ccal$ is a retract of such an $X$ and $\Dcal$ is idempotent-complete, we get that $\overline{f}:\Ccal\to\Dcal$, which concludes.
\end{proof}

\begin{rmq}
    Note that the above proof also shows that, if $\Ccal$ is generated by finite colimits by $G$ without needing retracts, then the idempotent-completeness is not necessary for $\Dcal$.
\end{rmq}

\begin{prop} \label{MapOfKernelsIsSplit}
    Let $F:\tB\Z\to\CatP$ be a functor valued in idempotent-complete categories. Then, the functor
    $$
        \begin{tikzcd}
            \Psi:\ker(\Phi^{\mathrm{p}}_{F\mid_{\TwAr(\Delta^1)}})\arrow[r] & \ker(\Phi^{\mathrm{p}}_{F})
        \end{tikzcd}
    $$ 
    is in the image of the functor $\Hyp:\CatEx\to\CatP$. Moreover, it admits a retraction $\Psi_r$ which is also in the image of $\Hyp$.
\end{prop}
\begin{proof}
    Both the source and the target of $\Psi$ are hyperbolic by respectively Proposition \ref{KernelHyperbolic} and Lemma \ref{KernelMetabolicCubes} for $n=1$. For the latter, we note the Lemma only proves metabolicity but it is clear that $\Lcal$ and $D_\Qoppa(\Lcal)$ are fully orthogonal so that the category is in fact hyperbolic.

    To check that the functor is in the image of $\Hyp$, it suffices to show that it sends the canonical Lagrangian of the hyperbolic category to the Lagrangian of the target (here, the canonical Lagrangian of $\Hyp(\Ccal)$ is none other $\Ccal$ and we are simply saying that $\Hyp(\Ccal)\to\Hyp(\Dcal)$ is of the form $\Hyp(f)$ if and only if it sends $\Ccal$ to $\Dcal$). 
    
    Write $\gamma_+$ and $\gamma_-$ for the two arrows of $\TwAr(\Delta^1)$ such that the former is mapped to $\alpha^+$ and the latter to $\alpha_-$ by $\TwAr(\Delta^1)\to\Ical$. We claim that $z^{\gamma^+}_x$ in the notation of Lemma \ref{KernelMetabolicCubes} is sent to $z^+_x$ in the notation of Proposition \ref{KernelHyperbolic}; this can be checked directly by inspecting the respective formulas for both sections. \\
    
    Furthermore, the induced functor
    $$
        \Psi(-):\Map(z^{\alpha^+}_x, z^{\alpha^+}_x)\longrightarrow\Map(z^+_x, z^+_x)
    $$
    has a retraction. This is because the right hand side is the space of squares
    $$
        \begin{tikzcd}
            X\arrow[d, "f"] &\arrow[l, "\can"] F(\alpha^-)(X)\arrow[d, "g"] && X\arrow[d, "f"] &\arrow[l, "0"] F(\alpha^-)(X)\arrow[d, "g"]\\
            X &\arrow[l, "\can"] F(\alpha^-)(X) && X &\arrow[l, "0"] F(\alpha^-)(X)
        \end{tikzcd}
    $$
    whereas the left hand side only has the left square as part of the datum. 
    
    The datum of a right hand square can equivalently be described as the datum of a map $\Sigma F(\alpha^-)(X)\to X$. Under this description, the functor $\Psi(-)$ sends the zero map $0:\Sigma F(\alpha^-)(X)\to X$ and consequently, it has a retraction given by forgetting this extra map. Since each Lagrangian is generated by one object, namely $z^{\alpha^+}_x$ resp. $z^+_x$, it follows from Lemma \ref{MapsFromStableCatGeneratedByOneElement} that this restriction extends to the whole Lagrangians, since both kernels are idempotent-complete by hypothesis on $F$ and Lemma \ref{IdempotentCompleteStabUn} After applying $\Hyp$, this provides the wanted retraction $\Psi_r$ to $\Psi$. 
\end{proof}

\begin{rmq}
    Note that this retraction is not compatible with the exact sequence of diagram \eqref{eq-ComparisonReducesToKernels}: it is clear that the composite
    $$
        \begin{tikzcd}
            \ker(\Phi^{\mathrm{p}}_{F})\arrow[r, "\Psi_r"] & \ker(\Phi^{\mathrm{p}}_{F\mid_{\TwAr(\Delta^1)}})\arrow[r] & \UnH(F\mid_{\Ical}) 
        \end{tikzcd}
    $$ 
    is not fully-faithful and therefore cannot coincide with the inclusion of the kernel of $\Phi^{\mathrm{p}}_{F}$.
\end{rmq}

\begin{defi}
    Let $F:\tB\Z\to\CatP$ be a functor. We let $NE^{\hyp}(F)$ be the following (split)-cofiber:
    $$
        NE^{\hyp}(F):=\Sigma \cofib(E(\Psi))
    $$ 
    The shift appears for cosmetic reasons.
\end{defi}

Since $\Psi$ is a map in the image of $\Hyp$, $NE^{\hyp}(F)$ only depends on $E^{\hyp}$, the hyperbolisation of $E$. In particular, it vanishes for all bordism-invariant $E$. Note also that $NE^{\hyp}(F)$ does not coincide with $\Sigma E(\coker(\Psi))$; in fact, it is straightforward to check that $\coker(\Psi)\simeq 0$ as this holds before applying $\Hyp$.

\begin{lmm}
    Let $E:\CatP\to\Ecal$ be a Poincaré-Verdier localizing invariant and $F:\tB\Z\to\CatP$. Then, the composite
    $$
        \begin{tikzcd}
            \Omega NE^{\hyp}(F)\arrow[r] & E(\ker(\Phi^{\mathrm{p}}_{F}))\arrow[r] & E(\UnP(F\mid_{\Ical}))
        \end{tikzcd}
    $$
    is null. As a consequence, the splitting of Proposition \ref{MapOfKernelsIsSplit} is compatible with the maps in diagram \eqref{eq-ComparisonReducesToKernels}.
\end{lmm}
\begin{proof}
    Since $\ker(\Phi^{\mathrm{p}}_{F})$ is hyperbolic, the second map factors through $E^{\hyp}(\StabUn(F\mid_{\Ical}))$ and it is sufficient to check that this composite vanishes. But now, we are in a purely stable situation and it suffices to prove the vanishing of
    $$
        \begin{tikzcd}
            \Omega NG(F)\arrow[r] & G(\Lcal)\arrow[r] & G(\StabUn(F\mid_{\Ical}))
        \end{tikzcd}
    $$
    where $\Lcal$ is such that $\Hyp(\Lcal)\simeq\ker(\Phi^{\mathrm{p}}_{F})$ and $\Omega NG$ is now the split-cokernel of the map which gives the map of Proposition \ref{MapOfKernelsIsSplit} after applying $\Hyp$. This second vanishing holds more generally than for $G:=E^{\hyp}$. Indeed, in Theorem 3.15 of \cite{KirsteinKremer} and more precisely in its proof, this map is identified with
    $$
        \begin{tikzcd}[cramped]
            \Omega NG_{\alpha}(\Ccal)\oplus \Omega NG_{\alpha}(\Ccal)\arrow[r] &[-1em] \Omega NG_{\alpha}(\Ccal)\oplus \Omega NG_{\alpha}(\Ccal)\oplus G(\Ccal)\oplus G(\Ccal)\arrow[r, "{0\oplus 0\oplus(\id,\id)\oplus(\alpha,\id)}"] &[5em] G(\Ccal)\oplus G(\Ccal)
        \end{tikzcd}.
    $$
    where $NG_{\alpha}(\Ccal)$ is the nilterm of $\Ccal:=F(*)$ under $G$ twisted by the action of $\alpha:=F(1)$. In particular, this map is $\cyclicGrp{2}$-equivariant under the action that flips the copies of the nilterms and flips the copies of $(E^{\hyp})(F(*))$, which is precisely the action induced from the duality. Therefore, the map is still null after taking homotopy orbits.
\end{proof}

Assembling all of our results, we have therefore proven the following result:

\begin{thm}[Twisted Shaneson splitting --- Verdier-localizing case] \label{TSSVerdierLoc}
    Let $F:\tB\Z\to\CatP$ be an idempotent-complete Poincaré category with $\Z$-action, and $E:\CatP\to\Ecal$ a Poincaré--Verdier localizing invariant. Then, the the assembly map of $E$ fits in a split-exact sequence
    $$
        \begin{tikzcd}[cramped]
            \colim\limits_{\tB\Z}E\circ F\arrow[r] & E(\colim\limits_{\tB\Z} F)\arrow[l, bend right=25, dotted]\arrow[r] & NE^{\hyp}(F)
        \end{tikzcd}.
    $$
    Moreover, $NE^{\hyp}(F)$ vanishes if $E$ is bordism-invariant.
\end{thm}
\begin{proof}
    Using diagram \eqref{eq-ComparisonReducesToKernels}, we see that the cokernel of the assembly map for $E$ is exactly the once-shifted of the cokernel of $E(\Psi)$, which we have denote $NE^{\hyp}(F)$. Moreover, the assembly itself is the cokernel of maps with compatible retractions, therefore admits a retraction itself.
\end{proof}

Since Poincaré-Karoubi localizing invariants are invariant under idempotent-completion, we also have:

\begin{cor}[Twisted Shaneson splitting] \label{TSSKaroubiLoc}
    Let $F:\tB\Z\to\CatP$ be a Poincaré category with $\Z$-action and $E:\CatP\to\Ecal$ a Poincaré-Karoubi localizing invariant. Then, the the assembly map of $E$ fits in a split-exact sequence
    $$
        \begin{tikzcd}
            \colim\limits_{\tB\Z}E\circ F\arrow[r] & E(\colim\limits_{\tB\Z} F)\arrow[l, bend right=25, dotted]\arrow[r] & NE^{\hyp}(F)
        \end{tikzcd}.
    $$
    Moreover, $NE^{\hyp}(F)$ vanishes if $E$ is bordism-invariant.
\end{cor}

\bibliographystyle{alpha}
\bibliography{bibliographie}

\newcommand{\etalchar}[1]{$^{#1}$}
\begin{thebibliography}{CDH{\etalchar{+}}23d}

\bibitem[BEN]{BartelsEfimovNikolaus}
Arthur Bartels, Alexander Efimov, and Thomas Nikolaus.
\newblock In preparation.

\bibitem[BGN18]{BarwickGlasmanNardin}
Clark Barwick, Saul Glasman, and Denis Nardin.
\newblock Dualizing cartesian and cocartesian fibrations.
\newblock {\em Theory and Applications of Categories}, 33(4):67--94, 2018.

\bibitem[BHS64]{BHS-original}
Hyman Bass, Alex Heller, and Richard~G. Swan.
\newblock {The Whitehead group of a polynomial extension}.
\newblock {\em Inst. Hauted \'Etudes Sci. Publ. Math.}, 22, 1964.

\bibitem[CDH{\etalchar{+}}23a]{HermKI}
Baptiste Calmès, Emanuele Dotto, Yonathan Harpaz, Fabian Hebestreit, Markus
  Land, Kristian Moi, Denis Nardin, Thomas Nikolaus, and Wolfgang Steimle.
\newblock Hermitian {K}-theory for stable $\infty$-categories {I}:
  {F}oundations.
\newblock {\em Selecta Mathematica}, 29(1):1 -- 269, 2023.

\bibitem[CDH{\etalchar{+}}23b]{HermKII}
Baptiste Calmès, Emanuele Dotto, Yonathan Harpaz, Fabian Hebestreit, Markus
  Land, Kristian Moi, Denis Nardin, Thomas Nikolaus, and Wolfgang Steimle.
\newblock Hermitian {K}-theory for stable $\infty$-categories {II}: {C}obordism
  categories and additivity.
\newblock \url{https://arxiv.org/abs/2009.07224}, 2023.

\bibitem[CDH{\etalchar{+}}23c]{HermKIII}
Baptiste Calmès, Emanuele Dotto, Yonathan Harpaz, Fabian Hebestreit, Markus
  Land, Kristian Moi, Denis Nardin, Thomas Nikolaus, and Wolfgang Steimle.
\newblock {Hermitian {K}-theory for stable $\infty$-categories {III}:
  Grothendieck--Witt groups of rings}.
\newblock \url{https://arxiv.org/abs/2009.07225}, 2023.

\bibitem[CDH{\etalchar{+}}23d]{HermKIV}
Baptiste Calmès, Emanuele Dotto, Yonathan Harpaz, Fabian Hebestreit, Markus
  Land, Kristian Moi, Denis Nardin, Thomas Nikolaus, and Wolfgang Steimle.
\newblock Hermitian {K}-theory for stable $\infty$-categories {IV}: {P}oincaré
  motives and {K}aroubi-{G}rothendieck-{W}itt groups, 2023.
\newblock In preparation.

\bibitem[CDW24]{ChristDycherhoffWaldeLaxAdd}
Merlin Christ, Tobias Dycherhoff, and Tashi Walde.
\newblock Lax {A}dditivity.
\newblock \url{https://arxiv.org/abs/2402.12251}, 2024.

\bibitem[Cis19]{CisinskiBook}
Denis-Charles Cisinski.
\newblock {\em Higher Categories and Homotopical Algebra}, volume 180 of {\em
  Cambridge studies in advanced mathematics}.
\newblock Cambridge University Press, 2019.

\bibitem[Efi24]{EfimovDualizable}
Alexander~I. Efimov.
\newblock K-theory and localizing invariants of large categories.
\newblock \url{https://arxiv.org/abs/2405.12169}, 2024.

\bibitem[GHN17]{GepnerHaugsengNikolaus}
David Gepner, Rune Haugseng, and Thomas Nikolaus.
\newblock Lax colimits and free fibrations in $\infty$-categories.
\newblock {\em Documenta Mathematica}, 22:1225 -- 1266, 2017.

\bibitem[Goo91]{Goodwillie2}
Tom Goodwillie.
\newblock Calculus. {II}. {A}nalytic functors.
\newblock {\em K-Theory 5}, 4:295 -- 332, 1991.

\bibitem[Gra76]{Grayson}
Daniel Grayson.
\newblock Higher algebraic {K}-theory. {II} (after {Daniel Quillen}).
\newblock {\em Algebraic {K}-theory (Proc. Conf., Northwestern Univ., Evanston,
  Ill., 1976)}, 1976.

\bibitem[HHLN23]{HaugsengHebestreitLinskensNuiten}
Rune Haugseng, Fabian Hebestreit, Sil Linskens, and Joost Nuiten.
\newblock Lax monoidal adjunctions, two-variable fibrations and the calculus of
  mates.
\newblock {\em Proceedings of the London Mathematical Society},
  127(4):889--957, 2023.

\bibitem[Hin16]{Hinich}
Vladimir Hinich.
\newblock Dwyer–{K}an localization revisited.
\newblock {\em Homology, Homotopy and Applications}, 18(1):27--48, 2016.

\bibitem[HL13]{HopkinsLurie}
Michael Hopkins and Jacob Lurie.
\newblock {Ambidexterity in $K(n)$-Local Stable Homotopy Theory}.
\newblock
  \url{https://people.math.harvard.edu/~lurie/papers/Ambidexterity.pdf}, 2013.

\bibitem[KK24]{KirsteinKremer}
Dominik Kirstein and Christian Kremer.
\newblock A twisted {Bass-Heller-Swan} decomposition for localizing invariants.
\newblock \url{https://arxiv.org/abs/2410.22877}, 2024.

\bibitem[LT19]{LandTammePullbacks}
Markus Land and Georg Tamme.
\newblock {On the $K$-theory of pullbacks}.
\newblock {\em Annals of Mathematics}, 190(3):877 -- 930, 2019.

\bibitem[LT23]{LandTammePushouts}
Markus Land and Georg Tamme.
\newblock {On the $K$-theory of pushouts}.
\newblock \url{https://arxiv.org/abs/2304.12812}, 2023.

\bibitem[Lur08]{HTT}
Jacob Lurie.
\newblock {\em Higher Topos Theory}.
\newblock Princeton University Press, 2008.

\bibitem[Lur11]{LurieNotesLtheory}
Jacob Lurie.
\newblock {Course notes ``Algebraic L-theory and surgery'''}.
\newblock \url{https://www.math.ias.edu/~lurie/287x.html}, 2011.

\bibitem[Lur17]{HA}
Jacob Lurie.
\newblock Higher algebra.
\newblock \url{http://people.math.harvard.edu/~lurie/papers/HA.pdf}, 2017.

\bibitem[Lur18]{Kerodon}
Jacob Lurie.
\newblock Kerodon.
\newblock \url{https://kerodon.net}, 2018.

\bibitem[Lü25]{LuckIsoConjBook}
Wolfgang Lück.
\newblock Isomorphism conjectures in {K}- and {L}-theory.
\newblock \url{https://him-lueck.uni-bonn.de/data/ic.pdf}, 2025.

\bibitem[NS17]{NikolausScholze}
Thomas Nikolaus and Peter Scholze.
\newblock On topological cyclic homology.
\newblock {\em Acta Mathematica}, 221, 2017.

\bibitem[Ran73]{Ranicki-ALTIII}
Andrew~A. Ranicki.
\newblock {Algebraic L-theory III. Twisted Laurent extensions}.
\newblock In {\em Algebraic K-theory III - Battelle Institute Conference 1972},
  volume 343 of {\em {Lecture Notes in Mathematics}}, Springer, Berlin, 1973.

\bibitem[Ran81]{Ranicki-ESATS}
Andrew~A. Ranicki.
\newblock {\em {Exact sequences in the algebraic theory of surgery}}.
\newblock Mathematical Notes, Princeton University Press, 1981.

\bibitem[Ran92]{Ranicki-blue}
Andrew~A. Ranicki.
\newblock {\em {Algebraic L-theory and topological manifolds}}.
\newblock Cambridge Tracts in Mathematics 102, Cambridge University Press,
  1992.

\bibitem[Rez23]{MO-Rezk}
Charles Rezk.
\newblock {Answer to ``When do Kan extensions preserve limits/colimits?''}.
\newblock \url{https://mathoverflow.net/q/450214}, 2023.

\bibitem[RW12]{RanickiWeiss}
Andrew~A. Ranicki and Michael Weiss.
\newblock {On The Algebraic $L$-theory of $\Delta$-sets}.
\newblock {\em Pure and Applied Mathematics Quarterly}, 8(2):423--449, 2012.

\bibitem[Sau23]{SaunierFundamental}
Victor Saunier.
\newblock The fundamental theorem of localizing invariants.
\newblock {\em Annals of K-theory}, 8(4):609--643, 2023.

\bibitem[Sha69]{Shaneson}
Julius Shaneson.
\newblock {Wall's surgery obstruction groups for $G \times Z$}.
\newblock {\em Annals of Mathematics}, 90(2):296--334, 1969.

\bibitem[SW25]{SaunierWinges}
Victor Saunier and Christoph Winges.
\newblock On exact categories and their stable envelopes.
\newblock \url{https://arxiv.org/abs/2502.03408}, 2025.

\bibitem[Yan21]{Yanovski}
Lior Yanovski.
\newblock The monadic tower for $\infty$-categories.
\newblock {\em Journal of Pure and Applied Algebra}, 226, 2021.

\end{thebibliography}

\end{document}